\documentclass{amsart}
\usepackage{amssymb}
\usepackage{color}
\usepackage{verbatim} %Allows use of \begin{comment} ... \end{comment}
\newtheorem{Thm}[equation]{Theorem}
\newtheorem{Pro}[equation]{Proposition}
\newtheorem{cor}[equation]{Corollary}
\newtheorem{lem}[equation]{Lemma}
\newtheorem{Example}[equation]{Example}
\newtheorem{con}[equation]{Convention}
\newtheorem{deft}[equation]{Definition}
\newtheorem{rem}[equation]{Remark}

\numberwithin{equation}{section}

\address{%
}

\def\a1s{a_1,\cdots, a_s}
\def\a{\alpha}

\def\aa{\mathcal A}

\def\fa{\mathfrak{a}}

\def\andd{\quad\hbox{and}\quad}

\def\fb{\frak{b}}
\def\b{\beta}

\def\bb{\mathcal{B}}

\def\bl4{B_{\ell\geq4}}
\def\cc{{\mathcal C}}

\def\bbbc{{\mathbb C}}

\def\d{\delta}
\def\D{\Delta}

\def\dd{\mathcal D}

\def\End{\hbox{End}}

\def\bbbf{\mathbb{F}}

\def\gg{{\mathcal G}}

\def\fg{\mathfrak{g}}

\def\heart{\hbox{\tiny$\heartsuit$}}

\def\hh{{\mathcal H}}

\def\fh{\mathfrak{h}}

\def\lam{\lambda}
\def\Lam{\Lambda}
\def\LL{\mathcal{L}}

\def\fl{\mathfrak{L}}
\def\LLc{\mathcal{L}_c}

\def\ep{\epsilon}
\def\fm{(\cdot,\cdot)}

\def\bbbq{\mathbb{Q}}

\def\bbbn{\mathbb{N}}

\def\bbbr{{\mathbb R}}

\def\1k{\frac{1}{k}}
\def\op{\oplus}
\def\ot{\otimes}

\def\la{\langle}
\def\ra{\rangle}

\def\sub{\subseteq}
\def\sg{\sigma}

\def\rcross{R^{\times}}

\def\pf{\noindent{\bf Proof. }}

\def\rre{R_{re}}
\def\rim{R_{im}}

\def\fp{\mathfrak{p}}

\def\ss{\mathcal{S}}

\def\u{{\mathcal U}}

\def\v{{\mathcal V}}

\def\w{{\mathcal W}}

\def\bbbz{{\mathbb Z}}
\def\ss{\mathcal{S}}
 
\def\1il{1\leq i\leq\ell}

\begin{document}
\centerline{\bf Root Graded Lie Superalgebras}
\vspace{0.7cm}
\centerline{Malihe Yousofzadeh\footnote{ma.yousofzadeh@sci.ui.ac.ir,\\ Department of Mathematics, University of Isfahan, Isfahan, Iran,
P.O.Box 81745-163, and School of Mathematics, Institute for Research in
Fundamental
Sciences (IPM), P.O. Box: 19395-5746, Tehran, Iran.\\
This research was in part
supported by a grant from IPM (No. 91170415) and partially carried out in IPM-Isfahan branch.
}}

\vspace{2cm}
\textbf{Abstract.} We define  root graded Lie superalgebras and study their  connection with  centerless cores of extended affine Lie superalgebras; our definition  generalizes the known notions of root graded Lie superalgebras.
\setcounter{section}{-1}         %<-------------------
\section{Introduction}
Motivated by a construction
appearing  in the classification of finite dimensional simple Lie
algebras containing  nonzero  toral subalgebras \cite{Se},  S. Berman and R. Moody \cite{BM} introduced the notion of
a  Lie algebra graded by an irreducible  reduced finite  root
system. This notion was generalized to Lie algebras graded by a locally finite  root
system and well studied through a  variety  of papers; recognition theorems for root graded Lie algebras are found in \cite{BM}, \cite{BZ},   \cite{N2}, \cite{ABG2}, \cite{BS}, \cite{you1} and their  central extensions  have been studied in \cite{ABG1}, \cite{ABG2} and \cite{you4}. Roughly speaking, a root graded Lie algebra is a Lie algebra which is  graded by  the root lattice of an irreducible locally finite root system $R$ and  contains a  locally finite split simple Lie algebra  whose root system is a  full subsystem of $R.$ One of the important phenomena in  the study of root graded Lie algebras is their interaction  with other classes of Lie algebras such as invariant affine reflection algebras \cite{N1}(see also \cite{AABGP}, \cite{AKY} and \cite{MY}); more precisely, the main ingredient in constructing  an invariant affine reflection algebra is a root graded Lie algebra \cite[\S 6]{N1}.

There
 have been two different approaches to define   root graded Lie superalgebras. One is working with Lie superalgebras which are  graded  by the root lattice of a locally finite root system and satisfy modified properties of a root graded Lie algebra \cite{N2}; the other one is working with a Lie superalgebra $\LL$ containing  a basic classical Lie superalgebra with a Cartan subalgebra $\hh$ with respect to which $\LL$ has a weight space decomposition satisfying certain properties \cite{BE1}. In fact in the latter  case, $\LL$ is graded by the root lattice of  the root system of a basic classical Lie superalgebra.  Root systems of  basic classical Lie superalgebras are exactly  generalized root systems introduced by  V. Serganova in 1996 \cite{serg}. Generalized root systems are called  finite root supersystems in \cite{you2} where  the author  introduces locally finite root supersystems and gives their classification. Locally finite root supersystems which  are extended by abelian groups  appear as the root systems of  specific Lie superalgebras named  extended affine Lie superalgebras \cite{you3}.    The so called core  of an extended affine Lie superalgebra is a Lie superalgebra  satisfying  certain properties which are in fact a super version of the features defining a root graded  Lie algebra. This motivates us to define root graded Lie superalgebras in a general setting. Our definition is a generalization of both mentioned   notions of root graded Lie superalgebras.
 In a series of papers, G. Benkart and A. Elduque  studied Lie superalgebras graded by finite root supersystems $C(n), D(m,n), D(2,1;\a), F(4), G(3),$ $A(m,n)$ and $B(m,n);$ see \cite{BE1}, \cite{BE2} and \cite{BE3}.
 We give a recognition theorem  for Lie superalgebras graded by the locally finite root supersystem of type $BC(I,J).$

This paper has been organized as follows. We begin the first section with gathering some information   regarding the locally finite Lie superalgebra ${\mathfrak{osp}}(I,J)$ and conclude the section with a separate subsection devoted to extended affine Lie superalgebras and their root systems. The material of this section are used to prove our recognition theorem for $BC(I,J)$-graded Lie superalgebras. In Section 2, we define root graded Lie superalgebras and realize  extended affine Lie superalgebras using root graded Lie superalgebras;  we  consider it as a first step of constructing  extended affine Lie superalgebras. Last section is exclusively devoted to the study of $BC(I,J)$-graded Lie superalgebras.

\section{Preliminaries}
Throughout this work, $\bbbf$ is a field of characteristic zero. Unless otherwise mentioned, all vector spaces are considered over $\bbbf.$ We denote the dual space of a vector space $V$ by $V^*.$ If $V$ is a vector space graded by an abelian group, we denote the degree of a homogeneous element $x\in V$ by $|x|;$ we also  make a convention that if for an element $x$ of $V,$  $|x|$  appears in an expression, by default, we assume that $x$ is homogeneous.  For a superspace $V,$ by $\hbox{End}_\bbbf(V),$ we mean the superspace of linear
endomorphisms of $V.$ If $A$ is an abelian group, we denote the group of automorphisms of $A$ by $Aut(A)$ and for a subset  $X$  of $A,$ by $\la X\ra,$ we mean the subgroup of $A$ generated by $X.$  Also we denote the cardinal number of a  set $S$ by $|S|;$ and for two symbols $i,j,$ by $\d_{i,j},$ we mean the Kronecker delta. We use $\uplus$ to indicate the disjoint union and  for a map $f:A\longrightarrow B$ and $C\sub A,$ by $f\mid_{_C},$ we mean the restriction of $f$ to $C.$
Finally, we denote the center of a Lie superalgebra $\gg$ by $Z(\gg)$ and for $\gg$-modules  $V$ and $W,$  by a $\gg$-module homomorphism from $V$ to $W,$ we mean  a  linear map $\varphi:V\longrightarrow W$ satisfying  $$\varphi(xv)=x\varphi(v);\;\; x\in\gg,\;\; v\in V.$$
%We note that if $V=V_0\op V_1$ is a $\gg$-module, then $W=W_0\op W_1$ in which $W_0:=V_1$ and $W_1=V_0$ is also a $\gg$-module isomorphic to $V.$
%Suppose that  $V$ and $W$ are two irreducible isomorphic $\gg$-module such that as $\gg_0$-modules $V_0\not\simeq V_1$, then there is a linear map $\phi:V\longrightarrow W$ with $\phi(xv)=x\phi(v)$ for all $x\in\gg$ and  $v\in  V.$ Then $\phi=\phi_0+\phi_1.$ But it is easy to see that $\phi_i$ is a $\gg$-module  homomorphism and so $\gg_0$-homomorphism. Since $ker\phi_i$ is a submodule of $V$ and $V$ is irreducible, we get that either $\phi_i$
%is a linear isomorphism or it is zero. If $\phi_i\neq 0$ for both $i=0$ and $i=1,$ we get that as $\gg_0$-modules, $V_0\simeq W_0,$ $V_1\simeq W_1$ (via $\phi_0$),  $V_0\simeq W_1,$ and $V_1\simeq W_0$ (via $\phi_1$). This is a contradiction, so at least one of $\phi_0$ and $\phi_1$ are zero and so $\phi$ is homogeneous.

\subsection{On locally finite Lie superalgebra ${\mathfrak{osp}}(I,J)$}\label{bc}
For two disjoint nonempty index sets  $I,J,$ suppose that $\{0,i,\bar i,\mid i\in I\cup J\}$ is a superset
with $|0|=|i|=|\bar i|=0$ for $i\in I$ and $|j|=|\bar j|=1$ for $j\in J.$  Take $\mathfrak{u}$ to be a vector superspace with a basis
$\{v_i\mid i\in I\cup\bar I\cup J\cup\bar J\cup\{0\}\}$  and $$|v_i|:=|i|;\;\;i\in I\cup\bar I\cup J\cup\bar J\cup\{0\}$$ in which by $\bar I$ (resp. $\bar J$), we mean $\{\bar i\mid i\in I\}$ (resp. $\{\bar j\mid j\in J\}$). Take $\fm$ to be the skew supersymmetric bilinear form $\fm$ on $\mathfrak{u}$ defined by   \begin{equation}\label{form1}(v_0,v_0):=1,(v_i,v_j):=0, (v_{\bar i},v_{\bar j}):=0,(v_i, v_{\bar j})=(-1)^{|i||j|}(v_{\bar j},v_i):=\d_{i,j}\end{equation} for $i,j\in I\cup J.$
Now for $j,k\in\{0\}\cup I\cup \bar I\cup J\cup \bar J,$
define
\begin{equation}\label{elementary2}e_{j,k}:\mathfrak{u}\longrightarrow\mathfrak{u};\;\;
v_i\mapsto \d_{k,i}v_j,\;\;\; (i\in \{0\}\cup I\cup \bar I\cup J\cup \bar J).\end{equation}  Then
\begin{equation}\label{elemnatry2}\mathfrak{gl}:=\hbox{span}_\bbbf\{e_{j,k}\mid j,k\in \{0\}\cup I\cup \bar I\cup J\cup \bar J\}\end{equation} is a
Lie subsuperalgebra of $\End_\bbbf(\mathfrak{u}).$
%Next set  $$\begin{array}{l}\gg_{\bar i}:=\{X\in \mathfrak{gl}_i\mid (Xv,w)=-(-1)^{|X||v|}(v, Xw);\forall v,w\in \u\};\;i=0,1\\
% \gg:=\gg_{\bar 0}\op\gg_{\bar 1}.\end{array}$$
%and
% $$\begin{array}{l}\ss_{\bar i} := \{s \in \mathfrak{gl}_i\mid  (su , v) = (-1)^{|s||u|}(u,sv),\; \forall u, v \in \u \andd str(s) = 0\};\;i=0,1\\
% \ss:=\ss_{\bar 0}\op\ss_{\bar 1}.\end{array}$$
Suppose that $\gamma\in\{1,-1\}.$ For $i=0,1,$ take  $$(\aa_\gamma)_{\bar i} := \{s \in \mathfrak{gl}_{\bar i}\mid str(s) = 0\andd  (su , v) = \gamma(-1)^{|s||u|}(u,sv),\; \forall u, v \in \mathfrak{u}\}$$ and  set
 $$\aa_\gamma:=(\aa_\gamma)_{\bar 0}\op(\aa_\gamma)_{\bar 1}.$$
 We next put \begin{equation}\label{gs}\mathfrak{g}:=\aa_{-1}\andd \mathfrak{s}:=\aa_1.\end{equation}
 One knows that $\fg$ is a Lie subsuperalgebra of $\End_\bbbf(\mathfrak{u}).$
Set\begin{equation}\label{h}\fh:=\hbox{span}_\bbbf\{h_t,d_k\mid t\in I,\;k\in J\}\end{equation} in which for  $t\in I$ and $ k\in J,$ $$h_t:=e_{t,t}-e_{\bar t,\bar t}\andd d_k:=e_{k,k}-e_{\bar k,\bar k}$$ and for $i\in I$ and $j\in J,$ define
$$\begin{array}{lll}
\begin{array}{c}
\ep_i:\fh\longrightarrow \bbbf\\
h_t\mapsto \d_{i,t},\;\;\; d_k\mapsto 0,
\end{array}&&
\begin{array}{c}
\d_j:\fh\longrightarrow \bbbf\\
h_t\mapsto 0,\;\;\; d_k\mapsto \d_{j,k},
\end{array}
\end{array}$$ in which $t\in I$ and $k\in J.$
Then  $\mathfrak{u}$ is a $\fg$-module equipped with a weight space decomposition $\mathfrak{u}=\op_{\a\in \D_\mathfrak{u}}\mathfrak{u}^\a$ with respect to $\fh,$ where \begin{equation}\label{delta}\D_\mathfrak{u}=\{0,\pm\ep_i,\pm\d_j\mid i\in I,j\in J\}\end{equation} with $$\mathfrak{u}^0=\bbbf v_0,\;\;\mathfrak{u}^{\ep_i}=\bbbf v_i,,\;\;\mathfrak{u}^{-\ep_i}=\bbbf v_{\bar i},\;\mathfrak{u}^{\d_j}=\bbbf v_j,,\;\;\mathfrak{u}^{-\d_j}=\bbbf v_{\bar j}$$
for $i\in I$ and $j\in J.$
Also  for $\gamma\in\{1,-1\},$ $\aa_\gamma$ is a $\fg$-module  having  a weight space decomposition with respect to $\fh.$
Taking $R$ (resp. $\D_{\mathfrak{s}}$) to be the set of weights of $\aa_{-1}$ (resp. $\aa_1$) with respect to $\fh,$  we have
\begin{equation}\label{root}\begin{array}{l}R=\{\pm\ep_r,\pm(\ep_r\pm\ep_s),\pm\d_p,\pm(\d_p\pm\d_q),\pm(\ep_r\pm\d_p)\mid  r, s\in I,\;p,q\in J,\;r\neq s\},\\
\\\D_{\mathfrak{s}}=\{\pm\ep_r,\pm(\ep_r\pm\ep_s),\pm\d_p,\pm(\d_p\pm\d_q),\pm(\ep_r\pm\d_p)\mid  r, s\in I,\;p,q\in J,\;p\neq q\}.\end{array}\end{equation}
Moreover, for $r, s\in I,\;p,q\in J,\;r\neq s$ and $p\neq q,$ we have
$$\begin{array}{ll}
(\aa_\gamma)^{\ep_r}=\hbox{span}_\bbbf  (e_{r,0}+\gamma e_{0,\bar r}),&
(\aa_\gamma)^{-\ep_r}=\hbox{span}_\bbbf  (e_{\bar r,0}+\gamma e_{0,r}),\\\\
(\aa_\gamma)^{\ep_r+\ep_s}=\hbox{span}_\bbbf  (e_{r,\bar s}+\gamma e_{s,\bar r}),&
(\aa_\gamma)^{-\ep_r-\ep_s}=\hbox{span}_\bbbf  (e_{\bar r,s}+\gamma e_{\bar s,r}),\\\\
(\aa_\gamma)^{\ep_r-\ep_s}=\hbox{span}_\bbbf  (e_{r,s}+\gamma e_{\bar s,\bar r)},&
(\aa_\gamma)^{2\ep_r}=\hbox{span}_\bbbf\d_{\gamma,1} e_{r,\bar r},\\\\
(\aa_\gamma)^{-2\ep_r}=\hbox{span}_\bbbf \d_{\gamma,1}e_{\bar r,r},&
(\aa_\gamma)^{2\d_p}=\hbox{span}_\bbbf\d_{\gamma,-1} e_{p,\bar p},\\\\
(\aa_\gamma)^{-2\d_p}=\hbox{span}_\bbbf \d_{\gamma,-1}e_{\bar p,p},&
(\aa_\gamma)^{\d_p+\d_q}=\hbox{span}_\bbbf (e_{p,\bar q}-\gamma e_{p,\bar q}),\\\\
(\aa_\gamma)^{-\d_p-\d_q}=\hbox{span}_\bbbf (e_{\bar p,q}-\gamma e_{\bar q,p}),&
(\aa_\gamma)^{\d_p-\d_q}=\hbox{span}_\bbbf (e_{p, q}+\gamma e_{\bar q,\bar p}),\\\\
(\aa_\gamma)^{\d_p}=\hbox{span}_\bbbf (e_{0,\bar p}-\gamma e_{p,0}),&
(\aa_\gamma)^{-\d_p}=\hbox{span}_\bbbf (e_{0,p}+\gamma e_{\bar p,0}),\\\\
(\aa_\gamma)^{\ep_r+\d_p}=\hbox{span}_\bbbf (e_{r,\bar p}+\gamma e_{p,\bar r}),&
(\aa_\gamma)^{-\ep_r-\d_p}=\hbox{span}_\bbbf (e_{\bar r,p}-\gamma e_{\bar p,r}),\\\\
(\aa_\gamma)^{\ep_r-\d_p}=\hbox{span}_\bbbf (e_{r,p}-\gamma e_{\bar p,\bar r}),&
(\aa_\gamma)^{-\ep_r+\d_p}=\hbox{span}_\bbbf (e_{\bar r,\bar p}+\gamma e_{p,r}).
\end{array}$$

\bigskip

In the literature,  $\fg$  is denoted by $\frak{osp}(I,J)$ and referred  to  as an {\it orthosymplectic} Lie superalgebra. We also refer to $\fg$  as  the  {\it split locally finite Lie superalgebra of type $B(I,J)$} ($B(m,n)$ if $|I|=m,|J|=n$) and say $\fh$ is the standard  {\it splitting Cartan subalgebra} of $\fg.$ We also refer to the $\fg$-module $\mathfrak{u}$ as the {\it natural module} of $\fg$ and to the $\fg$-module  $\mathfrak{s}$  as the {\it second natural module} of $\fg.$
We take $$\hbox{\small $\fg_{_B}:=\fg\cap \hbox{span}_\bbbf\{e_{i,j}\mid i,j\in I\cup \bar I\cup\{0\}\}\andd\fg_{_C}:=\fg\cap \hbox{span}_\bbbf\{e_{i,j}\mid i,j\in J\cup \bar J\}.$}$$ Then $\fg_{_B}$ (resp. $\fg_{_C}$) is a locally finite  split simple Lie algebra of type $B_I$ (resp. $C_J$) with splitting Cartan subalgebra $\hbox{span}_\bbbf\{h_i\mid i\in I\}$ (resp. $\hbox{span}_\bbbf\{d_j\mid j\in J\}$)   and corresponding root system $\{0,\pm\ep_i,\pm(\ep_i\pm\ep_j)\mid i,j\in I,\; i\neq j\}$ (resp. $\{0,\pm2\d_p,\pm(\d_p\pm\d_q)\mid p,q\in J,\; p\neq q\}$) \cite{NS}.
Moreover, $$\mathfrak{s}_{_B}:=\{x\in \mathfrak{s}\cap \hbox{span}_\bbbf\{e_{i,j}\mid i,j\in I\cup \bar I\cup\{0\}\}\mid tr(x)=0\}$$ is  a $\fg_{_B}$-module and   $$\mathfrak{s}_{_C}:=\{x\in \mathfrak{s}\cap \hbox{span}_\bbbf\{e_{i,j}\mid i,j\in J\cup \bar J\}\mid tr(x)=0\}$$ is  a $\fg_{_C}$-module.
 We finally note that if  $$|I|=m,\;|J|=n\andd  \mathfrak{I}:=\frac{1}{2m+1}\sum_{i\in \{0\}\cup I\cup \bar I}e_{ii}+\frac{1}{2n}\sum_{i\in J\cup \bar J}e_{ii},$$ then we have $$\mathfrak{s}_{\bar 0}=\mathfrak{s}_{_B}\op\mathfrak{s}_{_C}\op\bbbf \mathfrak{I}.$$
\medskip

We  have the following straightforward proposition.

\begin{Pro}\label{com}
Use the same notation as in the text and suppose that $I'\sub I,J'\sub J.$ Consider (\ref{root}) and take {\small $$\begin{array}{l}R':=R\cap\hbox{span}_\bbbz\{\ep_i,\d_p\mid i\in I',j\in J'\}\hbox{ and }
S:=\D_\frak{s}\cap\hbox{span}_\bbbz\{\ep_i,\d_p\mid i\in I',j\in J'\}.\end{array}$$} Set {\small$$\begin{array}{l}
{\displaystyle \gg:=\bigoplus_{\a\in R\setminus\{0\}}\fg^\a\op\sum_{\a\in R\setminus\{0\}}[\fg^\a,\fg^{-\a}]}\;\;\hbox{and}\;\;
{\displaystyle \ss:=\bigoplus_{\a\in S\setminus\{0\}}\mathfrak{s}^\a\op\sum_{\a\in S\setminus\{0\}}\fg^\a\cdot\mathfrak{s}^{-\a}},
\end{array}$$} then we have the following:

(i) $\gg$ is a Lie subsuperalgebra of $\fg$ isomorphic to $\mathfrak{osp}(I',J').$

(ii) Consider $\mathfrak{s}$ as a $\gg$-module, then $\ss$ is a $\gg$-submodule of $\mathfrak{s}$ isomorphic to the second natural module of $\gg.$

(iii) Suppose that  $V$ is a $\fg$-module isomorphic to $\fg$ and set {\small $$W:=\bigoplus_{\a\in R'\setminus\{0\}} V^\a\op\sum_{\a\in R'\setminus\{0\}}\fg^\a\cdot V^{-\a},$$} then $W$ is a $\gg$-module isomorphic to $\gg.$

(iv)   If $K$ is a $\fg$-module isomorphic to $\mathfrak{s},$ then {\small$$T:=\bigoplus_{\a\in S\setminus\{0\}}K^\a\op\sum_{\a\in S\setminus\{0\}}\fg^\a\cdot K^{-\a}$$} is a $\gg$-module isomorphic to the second natural module of $\ss.$

(v) Set $\Gamma_1:=\{0,\pm\ep_i,\pm\d_j\mid i\in I',j\in J'\}.$ If $U$ is a $\fg$-module isomorphic to $\mathfrak{u},$ then {\small $$M:=\bigoplus_{\a\in \Gamma_1\setminus\{0\}}U^\a\op\sum_{\a\in \Gamma_1\setminus\{0\}}\fg^\a\cdot U^{-\a}$$} is a $\gg$-module isomorphic to the natural module of $\mathfrak{osp}(I',J').$
\end{Pro}

\subsection{The Lie superalgebra $\mathfrak{osp(m,n)}$} In this  subsection, we suppose the field $\bbbf$ is algebraically closed and gather some facts regarding finite dimensional orthosympletic Lie superalgebras.
We keep the same notations as in the previous subsection and suppose $I=\{1,\ldots,m\}$ and $J=\{1,\ldots,n\}.$ We  denote the set of $\fg_{\bar 0}$-module homomorphisms from a $\fg_{\bar 0}$-module $X$ to  a $\fg_{\bar 0}$-module $Y$ by $\hom_{\fg_{\bar 0}}(X,Y).$

\begin{Pro}\label{hom}
Suppose that $\frac{2n}{2m+1}\not\in \bbbz,$ then
$$\begin{array}{ll}
\hom_{\fg_{\bar 0}}(\fg_{\bar 1}\ot\mathfrak{u}_{\bar 0},\mathfrak{s}_{\bar 0})=\{0\},&
\hom_{\fg_{\bar 0}}(\fg_{\bar 1}\ot\mathfrak{u}_{\bar 1},\mathfrak{s}_{\bar 0})=\{0\},\\
\hom_{\fg_{\bar 0}}(\fg_{\bar 1}\ot\mathfrak{u}_{\bar 0},\mathfrak{s}_{\bar 1})=\{0\},&
\hom_{\fg_{\bar 0}}(\fg_{\bar 1}\ot\mathfrak{u}_{\bar 1},\mathfrak{s}_{\bar 1})=\{0\},\\
\hom_{\fg_{\bar 0}}(\fg_{\bar 1}\ot\mathfrak{s}_{\bar 0},\mathfrak{u}_{\bar 0})=\{0\},&
\hom_{\fg_{\bar 0}}(\fg_{\bar 1}\ot\mathfrak{s}_{\bar 1},\mathfrak{u}_{\bar 0})=\{0\},\\
\hom_{\fg_{\bar 0}}(\fg_{\bar 1}\ot\mathfrak{s}_{\bar 0},\mathfrak{u}_{\bar 1})=\{0\},&
\hom_{\fg_{\bar 0}}(\fg_{\bar 1}\ot\mathfrak{s}_{\bar 1},\mathfrak{u}_{\bar 1})=\{0\},\\
\hom_{\fg_{\bar 0}}(\fg_{\bar 1}\ot\mathfrak{u}_{\bar 0},\fg_{\bar 0})=\{0\},&
\hom_{\fg_{\bar 0}}(\fg_{\bar 1}\ot\mathfrak{u}_{\bar 1},\fg_{\bar 0})=\{0\},\\
\hom_{\fg_{\bar 0}}(\fg_{\bar 1}\ot\mathfrak{u}_{\bar 0},\fg_{\bar 1})=\{0\},&
\hom_{\fg_{\bar 0}}(\fg_{\bar 1}\ot\mathfrak{u}_{\bar 1},\fg_{\bar 1})=\{0\},\\
\hom_{\fg_{\bar 0}}(\fg_{\bar 1}\ot\fg_{\bar 0},\mathfrak{u}_{\bar 0})=\{0\},&
\hom_{\fg_{\bar 0}}(\fg_{\bar 1}\ot\fg_{\bar 1},\mathfrak{u}_{\bar 0})=\{0\},\\
\hom_{\fg_{\bar 0}}(\fg_{\bar 1}\ot\fg_{\bar 0},\mathfrak{u}_{\bar 1})=\{0\},&
\hom_{\fg_{\bar 0}}(\fg_{\bar 1}\ot\fg_{\bar 1},\mathfrak{u}_{\bar 1})=\{0\}.\\
\end{array}$$
\end{Pro}
\pf We first note that $\fg_{\bar 1}$ is a $\fg_{\bar 0}$-module isomorphic to $\mathfrak{u}_{\bar 0}\ot\mathfrak{u}_{\bar 1}$ and fix the  base $$\{\ep_1-\ep_2,\ldots,\ep_{m-1}-\ep_m,\ep_m,\d_1-\d_2,\ldots,\d_{n-1}-\d_n,2\d_n\}$$ for the root system of $\fg_{\bar 0}.$ With respect to this base, we denote the finite dimensional irreducible  $\fg_{\bar 0}$-module  of highest weight  $\lam$  by $V(\lam)$ and recall that
\begin{equation}\label{fact}
\parbox{4in}{
for two finite dimensional irreducible  highest weight modules $V(\lam)$ and $V(\mu),$ $V(\lam)\ot V(\mu)$ is decomposed  into finite dimensional irreducible highest weight modules of highest weights of the form  $\mu + \lam'$ for some  $\lam'$ in the set of weights of $V(\lam);$}
\end{equation}
see \cite[Exercise 24.12]{Hum}.
We also recall that if $V$ is an irreducible $\fg_{_B}$-module and $W$ is an irreducible $\fg_{_C}$-module, then $V\ot W$ is an irreducible $\fg_{\bar 0}$-module.
\medskip

\noindent $\underline{\bf \hom_{\fg_{\bar 0}}(\fg_{\bar 1}\ot\mathfrak{u}_{\bar 0},\mathfrak{s}_{\bar 0})=\{0\}}:$ Suppose that  $\mathfrak{u}_{\bar 0}\ot\mathfrak{u}_{\bar 0}=\op_{i=1}^r V_i$ is the decomposition of the $\fg_{\bar 0}$-module $\mathfrak{u}_{\bar 0}\ot\mathfrak{u}_{\bar 0}$ into  finite dimensional irreducible highest weight $\fg_{\bar 0}$-modules. Now we have
\begin{eqnarray*}
\hom_{\fg_{\bar 0}}(\fg_{\bar 1}\ot\mathfrak{u}_{\bar 0},\mathfrak{s}_{\bar 0})\simeq \hom_{\fg_{\bar 0}}((\mathfrak{u}_{\bar 0}\ot\mathfrak{u}_{\bar 0})\ot\mathfrak{u}_{\bar 1},\mathfrak{s}_{\bar 0})
&\simeq& \hom_{\fg_{\bar 0}}(\op_{i=1}^r V_i\ot\mathfrak{u}_{\bar 1},\mathfrak{s}_{\bar 0})\\
&\simeq&\op_{i=1}^r  \hom_{\fg_{\bar 0}}(V_i\ot\mathfrak{u}_{\bar 1},\mathfrak{s}_{\bar 0})\\
&\simeq&\op_{i=1}^r  \hom_{\fg_{\bar 0}}(V_i\ot\mathfrak{u}_{\bar 1},\mathfrak{s}_{_C})\\
&\op& \op_{i=1}^r  \hom_{\fg_{\bar 0}}(V_i\ot\mathfrak{u}_{\bar 1},\mathfrak{s}_{_B})\\
&\op&\op_{i=1}^r  \hom_{\fg_{\bar 0}}(V_i\ot\mathfrak{u}_{\bar 1},\bbbf \mathfrak{I}).
\end{eqnarray*}

If $\hom_{\fg_{\bar 0}}(V_i\ot\mathfrak{u}_{\bar 1},\mathfrak{s}_{_C})\neq \{0\}$ for some $i,$ then there is a nonzero $\fg_{\bar 0}$-module  homomorphism $\varphi\in \hom_{\fg_{\bar 0}}(V_i\ot\mathfrak{u}_{\bar 1},\mathfrak{s}_{_C})\neq \{0\}.$ But $V_i\ot\mathfrak{u}_{\bar 1}$ and $\mathfrak{s}_{_C}$ are irreducible, so $\varphi$ is an isomorphism.  We note that the  set of weights of $\mathfrak{s}_{_C}$ as a $\fg_{\bar 0}$-module is $\{0,\pm(\d_p\pm\d_q)\mid 1\leq p\neq q\leq n\}$ while  the set of weights of $V_i\ot\mathfrak{u}_{\bar 1}$ is a subset of $\{\pm\ep_i\pm\d_p ,\pm\ep_i\pm\ep_j\pm\d_p\mid 1\leq i,j\leq m,1\leq p\leq n\}$ which is a contradiction.
Using the same argument as above, we get that $\hom_{\fg_{\bar 0}}(V_i\ot\mathfrak{u}_{\bar 1},\mathfrak{s}_{_B})=\{0\}$ for all $1\leq i\leq r.$ Also as $\dim(V_i\ot\mathfrak{u}_{\bar 1})>1,$ there is no isomorphism from $V_i\ot\mathfrak{u}_{\bar 1}$ to $\bbbf \mathfrak{I}$ and so $\hom_{\fg_{\bar 0}}(V_i\ot\mathfrak{u}_{\bar 1},\bbbf \mathfrak{I})=\{0\}.$

\medskip

\noindent$ \underline{\bf \hom_{\fg_{\bar 0}}(\fg_{\bar 1}\ot\mathfrak{u}_{\bar 1},\mathfrak{s}_{\bar 0})=\{0\}}:$ Consider the decomposition  $\mathfrak{u}_{\bar 1}\ot\mathfrak{u}_{\bar 1}=\op_{i=1}^s V_i$ of the $\fg_{_C}$-module $\mathfrak{u}_{\bar 1}\ot\mathfrak{u}_{\bar 1}$ into irreducible submodules. We have
\begin{eqnarray*}
\hom_{\fg_{\bar 0}}(\fg_{\bar 1}\ot\mathfrak{u}_{\bar 1},\mathfrak{s}_{\bar 0})\simeq \hom_{\fg_{\bar 0}}(\mathfrak{u}_{\bar 0}\ot(\mathfrak{u}_{\bar 1}\ot \mathfrak{u}_{\bar 1}),\mathfrak{s}_{\bar 0})
&\simeq& \hom_{\fg_{\bar 0}}(\op_{i=1}^s \mathfrak{u}_{\bar 0}\ot  V_i,\mathfrak{s}_{\bar 0})\\
&\simeq&\op_{i=1}^s  \hom_{\fg_{\bar 0}}(\mathfrak{u}_{\bar 0}\ot  V_i,\mathfrak{s}_{\bar 0})\\
&\simeq&\op_{i=1}^s \hom_{\fg_{\bar 0}}(\mathfrak{u}_{\bar 0}\ot  V_i,\mathfrak{s}_{_B})\\
&\op& \op_{i=1}^s  \hom_{\fg_{\bar 0}}(\mathfrak{u}_{\bar 0}\ot  V_i,\mathfrak{s}_{_C})\\
&\op & \op_{i=1}^s  \hom_{\fg_{\bar 0}}(\mathfrak{u}_{\bar 0}\ot  V_i,\bbbf\mathfrak{I}).
\end{eqnarray*}

As before, since $\dim(\mathfrak{u}_{\bar 0}\ot  V_i)>1,$ there is no $\fg_{\bar0}$-module isomorphism from $\mathfrak{u}_{\bar 0}\ot  V_i$ to $\bbbf\mathfrak{I},$ so $\hom_{\fg_{\bar 0}}(\mathfrak{u}_{\bar 0}\ot  V_i,\bbbf\mathfrak{I})=\{0\}.$ Also the set of weights of $\mathfrak{u}_{\bar 0}\ot  V_i$ nontrivially intersects $\{\pm\d_p\pm\d_q\pm\ep_j\mid 1\leq p, q\leq m,1\leq j\leq n\}$ while the set of weights of $\fg_{\bar 0}$-module $\mathfrak{s}_{_B}$ and the set of weights  of  $\fg_{\bar 0}$-module $\mathfrak{s}_{_C}$  are $\{\pm\ep_i\pm\ep_j\mid 1\leq i,j\leq m\}$ and $\{0,\pm(\d_p\pm\d_q)\mid 1\leq p<q\leq n\}$ respectively. Therefore $\hom_{\fg_{\bar 0}}(\mathfrak{u}_{\bar 0}\ot  V_i,\mathfrak{s}_{_B})=\{0\}$ and $\hom_{\fg_{\bar 0}}(\mathfrak{u}_{\bar 0}\ot  V_i,\mathfrak{s}_{_C})=\{0\},$ for all $1\leq i\leq s.$

\medskip

\noindent$ \underline{\bf \hom_{\fg_{\bar 0}}(\fg_{\bar 1}\ot\mathfrak{s}_{\bar 0},\mathfrak{u}_{\bar 0})=\{0\}}:$ Suppose that $\mathfrak{u}_{\bar 1}\ot\mathfrak{s}_{_C}=\op_{i=1}^tV(\eta_i)$ is the decomposition of the $\fg_{_C}$-module  $\mathfrak{u}_{\bar 1}\ot\mathfrak{s}_{_C}$ into  irreducible submodules and note that by (\ref{fact}), $\{\eta_i\mid 1\leq i\leq t\}\sub\{(\d_{ 1}+\d_2)\pm\d_p\mid 1\leq p\leq n\},$ so for $1\leq i\leq t,$ $\eta_i\neq 0.$ Therefore, $\dim(V(\eta_i))\neq 1$ which in turn  implies that  $\dim(\mathfrak{u}_{\bar 0}\ot V(\eta_i))\neq \dim(\mathfrak{u}_{\bar 0}).$  In particular, since $\mathfrak{u}_{\bar 0}\ot V(\eta_i) $ and $\mathfrak{u}_{\bar 0}$ are irreducible $\fg_{\bar 0}$-modules, we have \begin{equation}\label{hom=zero}\hom_{\fg_{\bar 0}}(\mathfrak{u}_{\bar 0}\ot V(\eta_i),\mathfrak{u}_{\bar 0})=\{0\};\;1\leq i\leq t.\end{equation}  Next suppose that  $\mathfrak{u}_{\bar 0}\ot \mathfrak{s}_{_B}=\op_{i=1}^sV(\theta_i)$ is the decomposition of  the $\fg_{_B}$-module $\mathfrak{u}_{\bar 0}\ot \mathfrak{s}_{_B}$ into  finite dimensional irreducible highest weight submodules, then by (\ref{fact}),  $\{\theta_i\mid 1\leq i\leq s\}\sub\{2\ep_{ 1},2\ep_{1}\pm\ep_j\mid 1\leq j\leq m\}.$ Therefore, for each $1\leq i\leq s,$ the set of weights of $V(\theta_i)\ot\mathfrak{u}_{\bar 1}$  nontrivially intersects  $\{2\ep_{ 1}\pm\d_p,2\ep_{1}\pm\ep_j\pm\d_p\mid 1\leq j\leq m,1\leq p\leq n\}.$ So $V(\theta_i)\ot\mathfrak{u}_{\bar 1}$ is not isomorphic to $\mathfrak{u}_{\bar 0}$ or $\mathfrak{u}_{\bar 1}$ as the set of weights of $\mathfrak{u}_{\bar 0}$  is $\{0,\pm\ep_j\mid 1\leq j\leq m\}$ and the set of weights of $\mathfrak{u}_{\bar 1}$  is $\{\pm\d_p\mid 1\leq p\leq n\};$ in particular, since $V(\theta_i)\ot\mathfrak{u}_{\bar 1},\mathfrak{u}_{\bar 0}$ and $\mathfrak{u}_{\bar 1} $ are irreducible $\fg_{\bar 0}$-module,  we have \begin{equation}\label{hom=zero2}\hom_{\fg_{\bar 0}}(V(\theta_i)\ot\mathfrak{u}_{\bar 1},\mathfrak{u}_{\bar 0})=\{0\}\andd \hom_{\fg_{\bar 0}}(V(\theta_i)\ot\mathfrak{u}_{\bar 1},\mathfrak{u}_{\bar 1})=\{0\}.\end{equation} We also note that  $\fg_{\bar 1}$ and $\mathfrak{u}_{\bar 0}$ as well as  $\fg_{\bar 1}$ and $\mathfrak{u}_{\bar 1}$ are non-isomorphic  irreducible $\fg_{\bar 0}$-modules, so  we get that \begin{equation}\label{hom=zero3}\hom_{\fg_{\bar 0}}(\fg_{\bar 1}\ot\bbbf\mathfrak{I},\mathfrak{u}_{\bar 0})=\{0\}\andd \hom_{\fg_{\bar 0}}(\fg_{\bar 1}\ot\bbbf\mathfrak{I},\mathfrak{u}_{\bar 1})=\{0\} .\end{equation}
Now using (\ref{hom=zero})-(\ref{hom=zero3}), we have
{\small
\begin{eqnarray*}
\hom_{\fg_{\bar 0}}(\fg_{\bar 1}\ot\mathfrak{s}_{\bar 0},\mathfrak{u}_{\bar 0})\hspace{-3mm}&\simeq&\hspace{-3mm} \hom_{\fg_{\bar 0}}(\fg_{\bar 1}\ot\mathfrak{s}_{_B},\mathfrak{u}_{\bar 0})\op \hom_{\fg_{\bar 0}}(\fg_{\bar 1}\ot\mathfrak{s}_{_C},\mathfrak{u}_{\bar 0})\op\hom_{\fg_{\bar 0}}(\fg_{\bar 1}\ot\bbbf\mathfrak{I},\mathfrak{u}_{\bar 0})\\
&\simeq&\hspace{-3mm} \hom_{\fg_{\bar 0}}((\mathfrak{u}_{\bar 0}\ot\mathfrak{s}_{_B})\ot\mathfrak{u}_{\bar 1},\mathfrak{u}_{\bar 0})\op \hom_{\fg_{\bar 0}}(\mathfrak{u}_{\bar 0}\ot(\mathfrak{u}_{\bar 1}\ot\mathfrak{s}_{_C}),\mathfrak{u}_{\bar 0})\\
&\simeq&\hspace{-3mm} \hom_{\fg_{\bar 0}}(\op_{i=1}^sV(\theta_i)\ot\mathfrak{u}_{\bar 1},\mathfrak{u}_{\bar 0})\op \hom_{\fg_{\bar 0}}(\mathfrak{u}_{\bar 0}\ot\op_{i=1}^tV(\eta_i),\mathfrak{u}_{\bar 0})\\
&\simeq&\hspace{-3mm}\op_{i=1}^s \hom_{\fg_{\bar 0}}(V(\theta_i)\ot\mathfrak{u}_{\bar 1},\mathfrak{u}_{\bar 0})\op \op_{i=1}^t\hom_{\fg_{\bar 0}}(\mathfrak{u}_{\bar 0}\ot V(\eta_i),\mathfrak{u}_{\bar 0})\\
&=&\{0\}.
\end{eqnarray*}
}

\noindent $ \underline{\bf \hom_{\fg_{\bar 0}}(\fg_{\bar 1}\ot\mathfrak{s}_{\bar 0},\mathfrak{u}_{\bar 1})=\{0\}}:$ For this, we first note that if $0\neq \varphi \in \hom_{\fg_{\bar 0}}(\mathfrak{u}_{\bar 0}\ot V(\eta_i),\mathfrak{u}_{\bar 1}),$ then $\varphi$ is an isomorphism and so $\dim(V(\eta_i))=2n/(2m+1)\not\in\bbbz,$ a contradiction. This together with (\ref{hom=zero2}) and (\ref{hom=zero3}) implies that
{\small\begin{eqnarray*}
\hom_{\fg_{\bar 0}}(\fg_{\bar 1}\ot\mathfrak{s}_{\bar 0},\mathfrak{u}_{\bar 1})\hspace{-3mm}&\simeq& \hspace{-3mm}\hom_{\fg_{\bar 0}}(\fg_{\bar 1}\ot\mathfrak{s}_{_B},\mathfrak{u}_{\bar 1})\op \hom_{\fg_{\bar 0}}(\fg_{\bar 1}\ot\mathfrak{s}_{_C},\mathfrak{u}_{\bar 1})\op\hom_{\fg_{\bar 0}}(\fg_{\bar 1}\ot\bbbf\mathfrak{I},\mathfrak{u}_{\bar 1})\\
&\simeq&\hspace{-3mm} \hom_{\fg_{\bar 0}}((\mathfrak{u}_{\bar 0}\ot\mathfrak{s}_{_B})\ot\mathfrak{u}_{\bar 1},\mathfrak{u}_{\bar 1})\op \hom_{\fg_{\bar 0}}(\mathfrak{u}_{\bar 0}\ot(\mathfrak{u}_{\bar 1}\ot\mathfrak{s}_{_C}),\mathfrak{u}_{\bar 1})\\
&\simeq& \hspace{-3mm}\hom_{\fg_{\bar 0}}(\op_{i=1}^sV(\theta_i)\ot\mathfrak{u}_{\bar 1},\mathfrak{u}_{\bar 1})\op \hom_{\fg_{\bar 0}}(\mathfrak{u}_{\bar 0}\ot\op_{i=1}^tV(\eta_i),\mathfrak{u}_{\bar 1})\\
&\simeq&\hspace{-3mm}\op_{i=1}^s \hom_{\fg_{\bar 0}}(V(\theta_i)\ot\mathfrak{u}_{\bar 1},\mathfrak{u}_{\bar 1})\op \op_{i=1}^t\hom_{\fg_{\bar 0}}(\mathfrak{u}_{\bar 0}\ot V(\eta_i),\mathfrak{u}_{\bar 1})\\
&\simeq&\hspace{-3mm}\op_{i=1}^t\hom_{\fg_{\bar 0}}(\mathfrak{u}_{\bar 0}\ot V(\eta_i),\mathfrak{u}_{\bar 1})=\{0\}.
\end{eqnarray*}}

These together with \cite[\S 3]{BE3} and  the fact  that  $\mathfrak{s}_{\bar 1}$ is a $\fg_{\bar 0}$-module isomorphic to the $\fg_{\bar 0}$-module  $\mathfrak{g}_{\bar1},$ completes the proof.\qed

\bigskip
Recall (\ref{root}) and suppose  that $|I|=m,|J|=n.$
One knows that $$\Pi:=\{\d_1-\d_{2},\d_2-\d_3,\ldots,\d_{n-1}-\d_n,\d_n-\ep_1,\ep_1-\ep_2,\ldots,\ep_{m-1}-\ep_m,\ep_m\}$$ is a fundamental system for the root system  $R$ of the finite dimensional  basic classical simple Lie superalgebra $\fg$ with respect to  the positive system $$\{\d_p\pm \d_q,2\d_p,\d_p,\d_p\pm\ep_i,\ep_i\pm\ep_j,\ep_i\mid 1\leq i<j\leq m,\;1\leq p<q\leq n\}.$$
Set $\rho:=(1/2)\sum_{\a\in R_0^+}\a-(1/2)\sum_{\a\in R_1^+}\a,$ where $R_0^+$ (resp. $R_1^+$) is  the set of positive even (resp. odd) roots, then
 we know from \cite[(2.2)]{K2} that
\begin{equation}\label{cas}
\parbox{3.3in}{the {\it Casimir element} $\Gamma$ of $\fg$ acts on the highest weight $\fg$-module of highest weight $\lam$ as $(\lam,\lam+2\rho)\hbox{id}$ }
  \end{equation}where `` ${\rm id}$ '' indicates the identity map. Moreover, we have
\begin{equation}\label{casimir}(\lam,\lam+2\rho)=\left\{\begin{array}{ll}
-2(n-m)& \hbox{if $\lam=\d_1$}\\
-2-4(n-m)& \hbox{if $\lam=2\d_1$}\\
2-4(n-m)& \hbox{if $\lam=\d_1+\d_{2}.$}
\end{array}\right.\end{equation}
Using \cite[Thm. 2.14]{ch-w} (see also \cite[Thm. 8]{K2}), we get that  the  only
nonzero elements of \begin{equation}\label{si}\Psi:=R\cup\{\pm2\ep_i\mid 1\leq i\leq m\}\end{equation} which can be the highest weight for a finite dimensional irreducible $\fg$-module are
$$\begin{array}{ll}
2\d_1,\;\; \d_1+\d_{2},\;\; \d_1 & \hbox{if } n \geq 2,\\
2\d_1,\;\; \d_1+\ep_1,\;\; \d_1& \hbox{if $n = 1.$ }\end{array}$$
One  knows that up to isomorphism, the only finite dimensional irreducible $\fg$-module whose highest weight is $2\d_1$ (resp. $\d_1$) is $\fg$ (resp. $\frak{u}$). Also  up to isomorphism, $\frak{s}$ is the only finite dimensional irreducible $\fg$-module whose highest weight is  $\d_1+\d_{2}$ if $n\neq 1$ and $\ep_1+\d_1$ if $n=1.$
The following lemma and its corollary are a slight generalization of a result of \cite[\S 3]{BE3}.

\begin{lem}\label{com-red}
Let $n\neq1$ and consider (\ref{si}). Suppose that  $X$ is a finite dimensional $\fg$-module equipped with a weight space decomposition with resect to $\fh$ whose set of weights is contained in $\Psi.$ Suppose that $Y$ is an irreducible  $\fg$-submodule of $X$ isomorphic to one of the $\fg$-modules $\fg,$ $\frak{u},$ $\frak{s}$ or  the trivial module such that $X/Y$ is also an irreducible  $\fg$-module  isomorphic  to one of the  above $\fg$-modules, then $X$ is completely reducible.
\end{lem}
\pf For $x\in X,$  we denote the image of $x$ in $X/Y$ under the canonical epimorphism $\bar{}:X\longrightarrow X/Y$ by $\bar x.$ Since $Y$ and $X/Y$ are finite dimensional irreducible  $\fg$-modules, they are highest weight modules. Suppose that $\lam$ and $\mu$ are the highest weights of $Y$ and $X/Y$ respectively. We first suppose that $(\lam,\lam+2\rho)\neq (\mu,\mu+2\rho).$   If $r$ is an eigenvalue of the action of the Casimir element $\Gamma$ on $X,$ then there is a nonzero $x\in X$ with $\Gamma x=rx,$ so $\Gamma\bar x=r\bar x.$ This means that  either $\bar x=0$ or $r=(\mu,\mu+2\rho)$ by (\ref{cas}). In the former case, $x\in Y$ and so $r=(\lam,\lam+2\rho).$ Therefore, the only eigenvalues for the action  of $\Gamma$ on $X$ are $(\lam,\lam+2\rho)$ and $(\mu,\mu+2\rho);$ in particular $X=X_\lam\op X_\mu$ in which $X_\lam$ and $X_\mu$ are the generalized eigenspaces corresponding to $(\lam,\lam+2\rho)$ and $(\mu,\mu+2\rho)$ respectively. Since  $\Gamma$ is  a $\fg$-module homomorphism, $X_\lam$ and $X_\mu$ are $\fg$-submodules of $X$ with $Y\sub X_\lam,$ therefore, we have $\frac{X}{Y}=\frac{X_\lam}{Y}\op\frac{X_\mu+Y}{Y}.$ But the only eigenvalue for the action of $\Gamma$ on $X/Y$ is $(\mu,\mu+2\rho),$ so $X_\lam/Y=\{0\};$ i.e., $X_\lam=Y$ is an irreducible $\fg$-module. This also implies that  $X_\mu\simeq X/Y$ is  an irreducible $\fg$-module. Therefore,  $X=X_\lam\op X_\mu$ is completely reducible. This completes the proof in this case. So from now till the end of the proof,   we assume $(\lam,\lam+2\rho)=(\mu,\mu+2\rho).$  By (\ref{casimir}), one of the  following cases can happen:
\begin{itemize}
\item $Y$ is isomorphic to $X/Y,$
\item one of $Y$ and $X/Y$ is the trivial module and the  other one is isomorphic to $\mathfrak{u},$
\item one of $Y$ and $X/Y$ is isomorphic to $\fg$  and the other one is isomorphic to $\frak{u},$
\item one of $Y$ and $X/Y$ is isomorphic to $\frak{s}$  and the other one is isomorphic to $\frak{u}.$
\end{itemize}
 Using the same argument as in \cite[\S 3]{BE3} together with  Proposition \ref{hom}, we get that in the first case, $X$ is completely reducible and that the last three cases result in a contradiction but   for the convenience of readers, we carry out the proof for one case. Suppose that $Y$ is isomorphic to $\frak{s}$ and $X/Y$ is isomorphic to $\frak{u},$ then by (\ref{casimir}), $-2(n-m)=2-4(n-m)$ and so $\frac{2n}{2m+1}\not\in \bbbz.$ Consider $X$ as a $\fg_{\bar 0}$-module, then $X_{\bar 0}$ as well as  $X_{\bar 1}$ are  completely reducible $\gg_{\bar 0}$-modules and so for $i=0,1,$ there is a $\fg_{\bar 0}$-submodule $Z_{\bar i}$ of $X_{\bar i}$ with $X_{\bar i}=Y_{\bar i}\op Z_{\bar i}.$ Set $Z:=Z_{\bar 0}\op Z_{\bar 1}$ which is a $\bbbz_2$-graded subspace of $X.$ Since $\fg$-module  $X/Y$ is isomorphic  to $\frak{u},$
% as a $\gg_{\bar 0}$-module either $(x/y)_{\bar 0}\simeq \frak{u}_{\bar 0}$ and $(X/Y)_{\bar 1}\simeq \frak{u}_{\bar 1}$ or  $(X/Y)_{\bar 0}\simeq \frak{u}_{\bar1}$ and $(X/Y)_{\bar 1}\simeq \frak{u}_{\bar 0}$
% and so
   $Z$ as a $\fg_{\bar 0}$-module is isomorphic to $\frak{u}.$ So $X=Y_{\bar 0}\op Y_{\bar 1}\op Z_{\bar 0}\op Z_{\bar 1}$ is  a decomposition of $X$ into $\fg_{\bar 0}$-modules with either  $Z_{\bar 0}\simeq \frak{u}_{\bar 0}$ and $Z_{\bar 1}\simeq \frak{u}_{\bar 1}$ or  $Z_{\bar 0}\simeq \frak{u}_{\bar1}$ and $Z_{\bar 1}\simeq \frak{u}_{\bar 0}.$
Since by Proposition \ref{hom}, $$\begin{array}{ll}\hom_{\fg_{\bar 0}}(\fg_{\bar 1}\ot\mathfrak{u}_{\bar 0},\mathfrak{s}_{\bar 0})=\{0\},&
\hom_{\fg_{\bar 0}}(\fg_{\bar 1}\ot\mathfrak{u}_{\bar 1},\mathfrak{s}_{\bar 0})=\{0\},\\
\hom_{\fg_{\bar 0}}(\fg_{\bar 1}\ot\mathfrak{u}_{\bar 0},\mathfrak{s}_{\bar 1})=\{0\},&
\hom_{\fg_{\bar 0}}(\fg_{\bar 1}\ot\mathfrak{u}_{\bar 1},\mathfrak{s}_{\bar 1})=\{0\},\end{array}$$ it follows that for $i=0,1,$ $\fg_{\bar 1} Z_{\bar i}\sub Z,$
and so $
\fg Z\sub Z.
$ This together with the fact that $Z$ is a $\bbbz_2$-graded subspace of $X$ implies that $Z$ is a $\fg$-submodule of $X.$ Also as $\fg$-module  $X/Y$ is isomorphic to $\frak{u},$ $Z$ as a $\fg$-module is isomorphic to $\frak{u}.$ Therefore, $X$ is completely reducible.
\qed

\begin{cor}
\label{cor}
Suppose that  $X$ is a finite dimensional  $\fg$-module equipped with a weight space decomposition with respect to $\fh.$ If the set of weights of $X$ is a subset of $\Psi,$ then $X$ is completely reducible such that its irreducible constituents  are isomorphic to one of $\fg$-modules $\fg,$ $\frak{s},$ $\frak{u}$ or  the trivial $\fg$-module.
\end{cor}
\pf One knows that $X$ has a composition series, say $\{0\}=X_0\sub X_1\sub X_2\sub \cdots\sub X_t=X.$ For each $1\leq i\leq t,$ $X_i$ is an $\fh$-submodule of  $X$ and so it inherits the weight space decomposition of $X$ with respect to $\fh.$ This implies that the set of weights of $X_i$ is contained in $\Psi$ and so the set of weights of the  irreducible $\fg$-module $X_i/X_{i-1}$ is contained in $\Psi.$ Therefore, $X_i/X_{i-1}$ is a finite dimensional irreducible $\fg$-module  whose  highest weight is an element of $\Psi;$ in particular,  it either  is isomorphic to one of $\fg$-modules $\fg,$ $\frak{s},$ $\frak{u}$ or is the trivial $\fg$-module. Now the result follows using Lemma  \ref{com-red}.
\qed

\subsection{Extended affine Lie superalgebras and their root systems}\label{extended} In this subsection, we recall the notions of extended affine Lie superalgebras and extended affine root supersystems from \cite{you3}. We  prove Lemma 2.28 which is essential for the study of root graded Lie superalgebras.
In the sequel, by a {\it symmetric form} on an additive abelian group $A,$ we mean a map $\fm: A\times A\longrightarrow \bbbf$ satisfying
\begin{itemize}
\item $(a,b)=(b,a)$ for all $a,b\in A,$
\item $(a+b,c)=(a,c)+(b,c)$ and $(a,b+c)=(a,b)+(a,c)$ for all $a,b,c\in A.$
\end{itemize}
In this case, we set  $A^0:=\{a\in A\mid(a,A)=\{0\}\}$ and call it the {\it radical} of the form $\fm.$ The form is called {\it nondegenerate} if $A^0=\{0\}.$
We note that if the form is nondegenerate, $A$ is torsion free and we can identify $A$ as a subset of $\bbbq\ot_\bbbz A.$ In the following, if an abelian group  $A$ is  equipped with a nondegenerate symmetric form, we consider $A$ as a subset of $\bbbq\ot_\bbbz A$ without further explanation.
Also if $A$ is a vector space over $\bbbf,$ bilinear forms are used in the usual sense.
\medskip

We call a triple   $(\LL,\hh,\fm)$  a {\it super-toral triple} if
\begin{itemize}
\item {\rm $\LL=\LL_{\bar 0}\op\LL_{\bar 1}$ is  a nonzero Lie superalgebra, $\hh$  is a nontrivial subalgebra of $\LL_{\bar 0}$ and  $\fm $ is an invariant  nondegenerate   even supersymmetric bilinear form $\fm$ on $\LL,$}
\item
{\rm $\LL$ has a  weight space decomposition $\LL=\op_{\a\in\hh^*}\LL^\a$ with respect to $\hh$ via the adjoint representation; we note that as $\LL_{\bar 0}$ as well as $\LL_{\bar 1}$ are $\hh$-submodules of $\LL,$ we have $\LL_{\bar 0}=\op_{\a\in \hh^*}(\LL_{\bar 0})^\a$ and $\LL_{\bar 1}=\op_{\a\in\hh^*}(\LL_{\bar 1})^\a$  with $(\LL_{\bar i})^\a:=\LL_{\bar i}\cap \LL^\a,$ $ i=0,1,$}

 \item {\rm the restriction of the form  $\fm$ to $\hh\times \hh$ is nondegenerate.}
\end{itemize}
We call  $R:=\{\a\in \hh^*\mid \LL^\a\neq\{0\}\},$ the {\it root system} of $\LL$ (with respect to $\hh$). Each element of $R$ is called  a {\it root.} We refer to elements of $R_{ 0}:=\{\a\in \hh^*\mid (\LL_{\bar 0})^\a\neq\{0\}\}$ (resp. $R_{ 1}:=\{\a\in \hh^*\mid (\LL_{\bar 1})^\a\neq\{0\}\}$) as {\it even roots} (resp. {\it odd roots}). We note that $R=R_{ 0}\cup R_{ 1}.$ Suppose that $(\LL,\hh,\fm)$ is a super-toral triple with corresponding root system $R$ and take  $\mathfrak{p}:\hh\longrightarrow \hh^*$ to be  the function  mapping   $h\in\hh$ to $(h,\cdot).$    Since the form is nondegenerate on $\hh,$ the map $\mathfrak{p}$ is one to one. So for each element $\a$ of the image $\hh^\frak{p}$ of $\hh$ under the map $\mathfrak{p},$  there is a unique $t_\a\in\hh$ representing $\a$ through the form $\fm.$ Now we can transfer the form on $\hh$ to a form on $\hh^\mathfrak{p},$ denoted again by $\fm,$  and defined by \begin{equation}\label{form}(\a,\b):=(t_\a,t_\b)\;\;\;(\a,\b\in \hh^\fp).\end{equation}   It is proved that if  for  $\a\in R_i\setminus\{0\}$ ($i\in\{0,1\}$), there are $x_\a\in(\LL_{\bar i})^\a$ and $x_{-\a}\in(\LL_{\bar i})^{-\a}$ such that  $0\neq[x_\a,x_{-\a}]\in\hh,$ then $\a$ is  an element  of  $\hh^\mathfrak{p}$ \cite[Lem. 2.4]{you3}.
\begin{deft}
{\rm A super-toral triple  $(\LL=\LL_{\bar 0}\op\LL_{\bar 1},\hh,\fm)$ (or $\LL$ if there is no confusion),  with root system $R=R_{ 0}\cup R_{ 1},$ is called an {\it extended affine Lie superalgebra} if}
\begin{itemize}{\rm
\item{(1)}  for each  $\a\in R_i\setminus\{0\}$ ($i\in\{0,1\}$), there are $x_\a\in(\LL_{\bar i})^\a$ and $x_{-\a}\in(\LL_{\bar i})^{-\a}$ such that  $0\neq[x_\a,x_{-\a}]\in\hh,$
    \item{(2)}  for each  $\a\in R$ with $(\a,\a)\neq 0$ and $x\in \LL^\a,$ $ad x:\LL\longrightarrow\LL$ is a locally nilpotent linear transformation.}
\end{itemize}
\end{deft}

Suppose that $(\LL,\hh,\fm)$ is an extended affine Lie superalgebra with root system $R.$ It is proved  that for $\a\in R_i$ ($i=0,1$) with $(\a,\a)\neq 0,$ there are $e_\a\in (\LL_{\bar i})^\a,$ $f_\a\in(\LL_{\bar i})^{-\a}$ such that
$(e_\a,f_\a,h_\a:=\frac{2t_\a}{(\a,\a)})$ is an $\mathfrak{sl}_2$-super-triple in the sense that $$[e_\a,f_\a]=h_\a,[h_\a,e_\a]=2e_\a,[h_\a,f_\a]=-2f_\a.$$  Moreover, the subsuperalgebra $\gg(\a)$ of $\gg$ generated by $\{e_\a,f_\a,h_\a\}$ is either isomorphic to $\mathfrak{sl}_2$ or to $\mathfrak{spo}(2,1);$ see \cite{you3}.

\begin{deft}
{\rm Suppose that $(\LL,\hh,\fm)$ is an extended affine Lie superalgebra with root system $R.$
The subsuperalgebra of $\LL$ generated by $\LL^\a$ for $\a\in \{\b\in R\mid (\b, R)\neq\{0\}\}$ is called the {\it core} of $\LL.$}
\end{deft}

\begin{Example}
{\rm A basic classical finite dimensional simple Lie superalgebra $\LL$ is  an  extended affine Lie superalgebra with $\LL=\LL_c.$}
\end{Example}

By \cite[Pro. 3.3]{you3}, the root system of an extended affine Lie superalgebra is an extended affine root supersystem  in the following sense.

\begin{deft}\label{iarr}
{\rm Suppose that $A$ is a nontrivial additive abelian group, $R$ is a subset of $A$ and  $\fm:A\times A\longrightarrow \bbbf$ is  a symmetric form. Set
$$\begin{array}{l}
R^0:=R\cap A^0,\\
\rcross:=R\setminus R^0,\\
\rcross_{re}:=\{\a\in R\mid (\a,\a)\neq0\},\;\;\;\rre:=\rcross_{re}\cup\{0\},\\
\rcross_{im}:=\{\a\in R\setminus R^0\mid (\a,\a)=0 \},\;\;\; \rim:=\rcross_{im}\cup\{0\}.
\end{array}$$
We say $(A,\fm,R)$ is an {\it extended affine root supersystem}  if the following hold:
$$\begin{array}{ll}
(S1)& \hbox{$0\in R,$ and $\hbox{span}_\bbbz(R)= A,$}\\\\
%(S2)&\hbox{the form restricted to $A_{re}:=\hbox{span}_\bbbz\rre$ is nondegenerate,}\\\\
(S2)& \hbox{$R=-R,$}\\\\
(S3)&\hbox{for $\a\in \rre^\times$ and $\b\in R,$ $2(\a,\b)/(\a,\a)\in\bbbz,$}\\\\
(S4)&\parbox{4.5in}{ ({\it root string property}) for $\a\in \rre^\times$ and $\b\in R,$  there are nonnegative  integers  $p,q$  with $2(\b,\a)/(\a,\a)=p-q$ such that \begin{center}$\{\b+k\a\mid k\in\bbbz\}\cap R=\{\b-p\a,\ldots,\b+q\a\},$\end{center}
} \\\\
(S5)&\parbox{4.5in}{for $\a\in \rim$ and $\b\in R$ with $(\a,\b)\neq 0,$
$\{\b-\a,\b+\a\}\cap R\neq \emptyset.$ }
%note that if $(\b,\a)\neq 0,$ then $\b+r\a\in R$ if and only if $\b+r\a\in R^\times.$
%there are positive integers $p,q$ with $p+q\geq1$ such that $\{\b+k\a\mid k\in \bbbz\}\cap R=\{\b+k\a\mid -p\leq k\leq q\}$ and if for some $k\in\bbbz,$ both  $\b+k\a$ and $\b-k\a$ are elements of $R,$ then $(\b+k\a)/2,(\b-k\a)/2\in R.$
\end{array}$$
If there is no confusion, for the sake of simplicity, we say   {\it $R$ is an extended affine root supersystem in $A.$}
An extended affine root supersystem $R$ is called {\it irreducible} if $\rcross$ cannot be written as a disjoint union of two nonempty orthogonal subsets. An extended affine root supersystem $(A,\fm,R)$ is called a {\it locally finite root supersystem} if the form $\fm$ is nondegenerate.}
\end{deft}
\begin{Example}
{\rm  Extended affine root systems \cite{AABGP} and invariant affine reflection systems \cite{N1} are examples of extended affine root supersystems. Also a generalized root system \cite{serg} is a  locally finite root supersystem. }
\end{Example}

\begin{deft}
{\rm
Suppose that $(A,\fm,R)$ is  a locally finite root supersystem.
\begin{itemize}
\item A subset $S$ of  $R$ is called a  {\it sub-supersystem} if the restriction of the form to $\la S\ra$ is nondegenerate, $0\in S,$ for $\a\in S\cap\rre^\times, \b\in S$ and $\gamma\in S\cap\rim$ with $(\b,\gamma)\neq 0,$ $r_\a(\b)\in S$ and  $\{\gamma-\b,\gamma+\b\}\cap
    S\neq\emptyset.$
\item A sub-supersystem $S$ of $R$ is called  {\it full} if $\hbox{span}_\bbbq S=\bbbq\ot_\bbbz A.$
\item If  $(A,\fm, R)$ is irreducible, $R$ is said to be  of {\it real type} if  $\hbox{span}_\bbbq R_{re}=\bbbq\ot_{\bbbz} A;$ otherwise, we say it is of {\it imaginary type.}
\item  If $\{R_i\mid i\in I\}$ is a class of  sub-supersystems of $R$ which are mutually   orthogonal with respect the form $\fm$ and $R\setminus\{0\}=\uplus_{i\in I}(R_i\setminus\{0\}),$  we say $R$ is {\it the direct sum} of $R_i$'s and write $R=\op_{i\in I}R_i.$
\item The locally finite  root supersystem $(A,\fm,R)$ is called a {\it locally finite root system} if $\rim=\{0\};$ see \cite{LN}.
\end{itemize}}
\end{deft}

We have the following straightforward lemma; see \cite[Lem. 3.20]{you2}:
\begin{lem}
If $\{(X_i,\fm_i,S_i)\mid i\in I\}$ is a class of   locally finite root supersystems, then for   $X:=\op_{i\in I}X_i$ and $\fm:=\op_{i\in I}\fm_i$   $(X,\fm,S:=\cup_{i\in I} S_i)$ is a locally finite root supersystem. Also each locally finite root supersystem is a direct   sum of irreducible sub-supersystems.
\end{lem}

\begin{deft}
{\rm (i)  Two irreducible  extended affine  root supersystems $(A,\fm_1,R)$ and $(B,\fm_2,S)$ are called {\it isomorphic} if there is a group isomorphism $\varphi:A\longrightarrow B$ and a nonzero scalar $r\in\bbbf$ such that $\varphi(R)=S$ and  $(
a_1,a_2)_1=r(\varphi(a_1),\varphi(a_2))_2$ for all $a_1,a_2\in A.$

(ii) Suppose that  $(A,\fm,R)$  is an   extended affine root supersystem.  The subgroup $\w$ of $Aut(A)$
generated by $r_\a,$ $\a\in \rre^\times,$ is called the {\it Weyl group} of $R;$ we note that for $\a\in \rre^\times$ and  $a\in A,$ (S1) and (S3) imply that  $2(a,\a)/(\a,\a)\in\bbbz$ and so $r_\a:A\longrightarrow A$ mapping $a\in A$ to $a-\frac{2(a,\a)}{(\a,\a)}\a$ is a group  automorphism.
}
\end{deft}

\begin{Thm}(see \cite[\S4.14,
\S8]{LN} and \cite[Lem. 3.21]{you2})
Suppose that $T$ is a nonempty  index set  and $\u:=\op_{i\in
T}\bbbz\ep_i$ is the free $\bbbz$-module over   the
set $T.$ Define the symmetric  form $$\begin{array}{c}\fm:\u\times\u\longrightarrow\bbbf;\;\;
(\ep_i,\ep_j)=\d_{i,j}, \hbox{ for } i,j\in T
\end{array}$$
and set
\begin{equation}\label{locally-finite}
\begin{array}{l}
\dot A_T:=\{\ep_i-\ep_j\mid i,j\in T\}\hbox{ if $|T|>1$},\\
D_T:=\{\pm(\ep_i\pm\ep_j)\mid i,j\in T,\;i\neq j\} \hbox{ if $|T|>2$},\\
B_T:=\{\pm\ep_i,\pm(\ep_i\pm\ep_j)\mid  i,j\in T,\;i\neq j\},\\
C_T:=\{\pm2\ep_i,\pm(\ep_i\pm\ep_j)\mid  i,j\in T,\;i\neq j\},\\
BC_T:=B_T\cup C_T.
\end{array}
\end{equation}
These are irreducible locally finite root systems
in their $\bbbz$-span's. Moreover, each irreducible  locally finite root system is either  an irreducible finite root system  or an infinite locally finite root system  isomorphic to one of these locally finite root
systems.\end{Thm}

We refer to locally finite root systems listed in (\ref{locally-finite}) as  {\it type} $A,D,B,C$
and $BC$  respectively. We note that if $R$ is  an irreducible locally finite
root system as above, then  $(\a,\a)\in\bbbn$ for all
$\a\in R.$ This allows us to  define
$$\begin{array}{l}
R_{sh}:=\{\a\in R^\times\mid (\a,\a)\leq(\b,\b);\;\;\hbox{for all $\b\in R$} \},\\
R_{ex}:=R\cap2 R_{sh}\andd
R_{lg}:= R^\times\setminus( R_{sh}\cup R_{ex})\\
R_{red}:=\{0\}\cup R_{sh}\cup R_{lg}.
\end{array}$$
The elements of $R_{sh}$ (resp. $R_{lg},R_{ex}$ and $ R_{red}$) are called {\it
short roots} (resp. {\it long roots, extra-long roots} and {\it reduced roots}) of $R$. We point  out that following the usual notation in the literature,  the locally finite root system  of type $A$  is denoted by  $\dot A$ instead of $A,$ as all locally finite root systems listed above are spanning sets for $\bbbf\ot_\bbbz \u$ other than the one of type $A$ which spans a subspace of codimension 1.
\begin{con} {\rm We  make a convention that if a locally finite root system  $R$ is the direct sum of  subsystems $R_i,$ where $i$ runs over  a nonempty index set $I,$ for $*\in\{sh,lg,ex,red\},$ by $R_*,$ we mean $\cup_{i\in I}(R_i)_{*}.$}
\end{con}

\begin{Thm}[{\cite[Thm. 4.28]{you2}}]\label{classification I}
Suppose that $T,T'$ are index sets of cardinal numbers greater than $1$ with  $|T|\neq|T'|$ if $T,T'$ are both finite. Fix a symbol $\a^*$ and  pick $t_0\in T$ and $p_0\in T'.$ Consider the free $\bbbz$-module $X:=\bbbz\a^*\op\op_{t\in T}\bbbz\ep_t\op\op_{p\in T'}\bbbz\d_p$  and define the symmetric form $$\fm:X\times X\longrightarrow \bbbf$$ by $$\begin{array}{ll}
(\a^*,\a^*):=0,(\a^*,\ep_{t_0}):=1,(\a^*,\d_{p_0}):=1\\
(\a^*,\ep_t):=0,(\a^*,\d_q):=0&t\in T\setminus\{t_0\},q\in T'\setminus\{p_0\}\\
(\ep_t,\d_p):=0,(\ep_t,\ep_s):=\d_{t,s},(\d_p,\d_q):=-\d_{p,q}&t,s\in T,p,q\in T'.
\end{array}$$
 Take $R$ to be $\rre\cup \rim^\times$ as in the following table:
$${\small
\begin{tabular}{|c|c|c|}
\hline
type &$\rre$&$\rim^\times$
\\\hline
$\dot A(0,T)$& $\{\ep_t-\ep_s\mid t,s\in T\}$&$\pm\w\a^*$\\
\hline
$\dot C(0,T)$& $\{\pm(\ep_t\pm\ep_s)\mid t,s\in T\}$&$\pm\w\a^*$\\
\hline
$\dot A(T,T')$& $\{\ep_t-\ep_s,\d_p-\d_q\mid t,s\in T,p,q\in T'\}$&$\pm\w\a^*$\\
\hline
\end{tabular}}$$
in which  $\w$ is the subgroup of $Aut(X)$ generated by the reflections $r_\a$ $(\a\in \rre\setminus\{0\})$ mapping $\b\in X$ to $\b-\frac{2(\b,\a)}{(\a,\a)}\a,$ then $(A:=\la R\ra,\fm\mid_{A\times A}, R)$  is an irreducible locally finite root supersystem of imaginary type and conversely, each irreducible locally finite root supersystem of imaginary type  is isomorphic to one and only one of these root supersystems.
\end{Thm}

\begin{Thm}[{\cite[Thm. 4.37]{you2}}]
\label{classification II}
 Suppose  $(X_1,\fm_1,S_1),$ $\ldots,(X_n,\fm_n,S_n)$ for some $ n\in\{2, 3\},$  are irreducible locally finite root systems. Set $X:=X_1\op\cdots\op X_n$ and $\fm:=\fm_1\op\cdots\op\fm_n$ and consider  the locally finite root system  $(X,\fm,S:=\cup_{i=1}^n S_i).$ Take $\w$ to be the Weyl group of $S.$   If $1\leq i\leq n$ and $S_i$ is a finite root system of rank $\ell\geq 2,$ we take $\{\omega_1^i,\ldots,\omega_\ell^i\}\sub\bbbq\ot_\bbbz X_i$ to be a  set of fundamental weights for $S_i$ and if $S_i$ is one of infinite locally finite root systems $B_T, C_T, D_T$ or $BC_T$  as in (\ref{locally-finite}), by $\omega_1^i,$ we mean $\ep_1,$ where $1$ is a distinguished element of $T.$ Also if $S_i$ is  one of the  finite root systems $\{0,\pm\a\}$ of type $A_1$ or  $\{0,\pm\a,\pm2\a\}$ of type $BC_1,$ we set $\omega_1^i:=\frac{1}{2}\a.$
Consider $\d^*$ and $\dot R:=\dot R_{re}\cup\dot R_{im}^\times$ as in the following table:
$${\tiny
\begin{tabular}{|c|l|c|c|c|c|}
\hline
$ n$& $ S_i\;(1\leq i\leq n)$&$\dot R_{re}$&$\d^*$&$\dot R_{im}^\times$&\hbox{type}\\
\hline
$ 2$& $ S_1=A_\ell,\; S_2=A_\ell$ $(\ell\in\bbbz^{\geq1})$&$S_1\op S_2$&$\omega_1^1+\omega_1^2$&$\pm\w\d^*$&$A(\ell,\ell)$\\
\hline
$2$& $ S_1=B_T,\;S_2=BC_{T'}$ $(|T|,|T'|\geq2)$&$S_1\op S_2$&$\omega_1^1+\omega_1^2$&$\w\d^*$&$B(T,T')$\\
\hline
$2$& $ S_1=BC_T,\;S_2=BC_{T'}$ $(|T|,|T'|>1)$&$S_1\op S_2$&$\omega_1^1+\omega_1^2$&$\w\d^*$&$BC(T,T')$\\
\hline
$2$& $ S_1=BC_T,\;S_2=BC_{T'}$ $(|T|=1,|T'|=1)$&$S_1\op S_2$&$2\omega_1^1+2\omega_1^2$&$\w\d^*$&$BC(T,T')$\\
\hline
$2$& $ S_1=BC_T,\;S_2=BC_{T'}$ $(|T|=1,|T'|>1)$&$S_1\op S_2$&$2\omega_1^1+\omega_1^2$&$\w\d^*$&$BC(T,T')$\\
\hline
$2$& $ S_1=D_T,\;S_2=C_{T'}$ $(|T|\geq3,|T'|\geq2)$&$S_1\op S_2$&$\omega_1^1+\omega_1^2$&$\w\d^*$&$D(T,T')$\\
\hline
$2$& $ S_1=C_T,\;S_2=C_{T'}$ $(|T|,|T'|\geq2)$&$S_1\op S_2$&$\omega_1^1+\omega_1^2$&$\w\d^*$&$C(T,T')$\\
\hline
$2$& $ S_1=A_1,\;S_2=BC_{T}$ $(|T|=1)$&$S_1\op S_2$&$2\omega_1^1+2\omega_1^2$&$\w\d^*$&$B(1,T)$\\
\hline
$2$& $ S_1=A_1,\;S_2=BC_{T}$ $(|T|>1)$&$S_1\op S_2$&$2\omega_1^1+\omega_1^2$&$\w\d^*$&$B(1,T)$\\
\hline
$2$& $ S_1=A_1,\;S_2=C_T$ $(|T|\geq 2)$&$S_1\op S_2$&$\omega_1^1+\omega_1^2$&$\w\d^*$&$C(1,T)$\\
\hline
$2$& $ S_1=A_1,\;S_2=B_3$ &$S_1\op S_2$&$\omega_1^1+\omega_3^2$&$\w\d^*$&$AB(1,3)$\\
\hline
$2$& $ S_1=A_1,\;S_2=D_{T}$ $(|T|\geq3)$&$S_1\op S_2$&$\omega_1^1+\omega_1^2$&$\w\d^*$&$D(1,T)$\\
\hline
$2$& $ S_1=BC_1,\;S_2=B_T$ $(|T|\geq2)$&$S_1\op S_2$&$2\omega_1^1+\omega_1^2$&$\w\d^*$&$B(T,1)$\\
\hline
$2$& $ S_1=BC_1,\;S_2=G_2$ &$S_1\op S_2$&$2\omega_1^1+\omega_1^2$&$\w\d^*$&$G(1,2)$\\
\hline
$3$& $ S_1=A_1,\;S_2=A_1,\; S_3=A_1$ &$S_1\op S_2\op S_3$&$\omega_1^1+\omega_1^2+\omega_1^3$&$\w\d^*$&$D(2,1,\lam)$\\
\hline
$3$& $ S_1=A_1,\; S_2= A_1,\; S_3:=C_T$ $(|T|\geq2)$&$S_1\op S_2\op S_3$&$\omega_1^1+\omega_1^2+\omega_1^3$&$\w\d^*$&$D(2,T)$\\
\hline
\end{tabular}
}$$
For $1\leq i\leq n,$ normalize the form $\fm_i$ on $X_i$ such that  $(\d^*,\d^*)=0$ and that for type $D(2,T),$ $(\omega_1^1,\omega_1^1)_1=(\omega_1^2,\omega_1^2)_2.$  Then $(\la \dot R\ra,\fm\mid_{\dot X\times \dot X},\dot R)$ is an irreducible locally finite root supersystem of real type and conversely, if  $(\dot X,\fm,\dot R)$ is an irreducible locally finite root supersystem of real type, it is either an irreducible locally finite root system or  isomorphic to one and only  one of the locally finite root supersystems listed in the above table.
\end{Thm}

\begin{lem}\label{lem}
Suppose that $(\LL,\hh,\fm)$ is an extended affine Lie superalgebra with irreducible  root system $R.$ Set $\v:=\hbox{span}_\bbbf R$ and denote the induced form on $\v$ again by $\fm;$ see (\ref{form}). Take $\v^0$ to be the radical of the form. Suppose that  $\;\bar{}:\v\longrightarrow \bar \v:=\v/\v^0$ is  the canonical projection map and take  $\bar R$ to be  the image of $R$ under the projection map ``\;\;$\bar{}$\;\;''. Denote by $(\cdot,\cdot\bar ),$ the induced form on $\bar \v,$ then we have the following:

(i)  $(\bar A:=\la\bar R\ra,(\cdot,\cdot\bar)|_{\bar A\times\bar A},\bar R)$ is an irreducible  locally finite root supersystem.

(ii) There is a triple $(\dot \v,\dot R,\{S_{\dot\a}\}_{\dot\a\in\dot R})$  in which \begin{itemize}
\item $\dot \v$ is a subspace  of $\v$ with  $\v=\dot\v\op\v^0,$
\item $\dot R\sub \dot\v$ and $\dot R$ is a locally finite root supersystem (in its $\bbbz$-span) isomorphic to $\bar R;$ in particular, $\dot R_{re}$ is a locally finite root system,
\item for each $\dot\a\in \dot R,$ $S_{\dot \a}$ is a nonempty subset of $\v^0$ such that
\begin{itemize}
\item $R=\cup_{\dot\a\in\dot R}(\dot\a+S_{\dot\a}), $
\item  $0\in S_{\dot\a} \hbox{ for }
\dot\a\in \left\{\begin{array}{ll}
(\dot R_{re})_{red} & \hbox{\small $\dot R$ is of real type,}\\
\dot R & \hbox{\small $\dot R$ is of imaginary  type,}
\end{array}\right.$
\item  if $\dot R_{im}\neq\{0\}$ and $\dot R$ is of type $X\neq A(\ell,\ell),C(T,T'),C(1,T),$ then for all $\dot\a,\dot\b\in (\dot R_{re})_{sh},$   $S_{\dot\a}=S_{\dot \b};$ also for all $\dot\a,\dot\b\in (\dot R_{re})_{lg}\cup \dot R_{im}^\times,$  $S_{\dot\a}=S_{\dot \b},$
\item   if $\dot R_{im}\neq\{0\}$ and $\dot R$ is of type $X\neq A(\ell,\ell),C(T,T'),C(1,T),$ setting $S:=S_{\dot\a}$ for some $\dot\a\in (\dot R_{re})_{sh}$  and $F:=S_{\dot\b}$ for some  $\dot\b\in \dot R_{im},$ we get that $F$ is a subgroup of $\v^0$ and   $$S-2S\sub S,\;\;S+F\sub S\andd 2S+F\sub F.$$
\end{itemize}
\end{itemize}
\end{lem}
\pf Using the same argument as in \cite[Lem. 3.10]{you2}, one can see that $\bar R_{re}$ is locally finite in its $\bbbf$-span in the sense that it intersects each finite dimensional subspace of $\hbox{span}_\bbbf \bar R_{re}$ in a finite set. So using Lemmas 3.10, 3.12 and 3.21  of \cite{you2}, we get that $\bar R$ is an irreducible locally finite root supersystem in its $\bbbz$-span. Also using \cite[Lem. 3.5]{you2}; we get that $\bar R_{re}$ is a locally finite root system and the restriction of the form  $(\cdot,\cdot\bar)$ to $\bar\v_{re}:=\hbox{span}_\bbbf \bar R_{re}$ is nondegenerate. Therefore we have
\begin{equation}\label{non}
\parbox{4.3in}{the restriction of the form $(\cdot,\cdot\bar)$ to $\bar\v_{\bbbq}:=\hbox{span}_\bbbq \bar R_{re}$ is nondegenerate.}
 \end{equation} Since  $\bar R_{re}$ is a locally finite root system, by \cite[Lem. 5.1]{LN2},  it contains a $\bbbz$-linearly independent subset $T$ such that \begin{equation}\label{ref}\w_TT=(\bar R_{re})^\times_{red}=\bar R_{re}\setminus\{2\bar \a\mid\a\in R_{re}\},\end{equation} in which by $\w_T,$ we mean  the subgroup of the Weyl group of $\bar R_{re}$ generated by $r_{\bar\a}$ for all $\bar\a\in T.$  On the other hand, we know there is a subset $\Pi$ of $R$ such that $\bar \Pi$ is  a $\bbbz$-basis for $\hbox{span}_\bbbz\bar R;$ see \cite[Lem. 3.13]{you3}. This allows us to define the linear isomorphism  $$\varphi:\hbox{span}_\bbbq \bar R\longrightarrow \bbbq\ot_\bbbz \hbox{span}_\bbbz\bar R$$ mapping $\bar \a$ to $1\ot\bar \a$ for all $\a\in\Pi.$
Now suppose that  $\bar R$ is of real type, then
\begin{eqnarray*}
\varphi(\hbox{span}_\bbbq \bar R_{re})
=
\hbox{span}_\bbbq (1\ot \bar R_{re})=\bbbq\ot\hbox{span}_\bbbz\bar R=\varphi(\hbox{span}_\bbbq \bar R)
\end{eqnarray*}
which in turn implies that $\hbox{span}_\bbbq \bar R=\hbox{span}_\bbbq \bar R_{re}.$ Therefore,  $\hbox{span}_\bbbq \bar R=\hbox{span}_\bbbq T$ and so $\hbox{span}_\bbbf \bar R=\hbox{span}_\bbbf T .$ But $T$ is $\bbbz$-linearly independent and so it is $\bbbq$-linearly independent. We now prove that $T$ is $\bbbf$-linearly independent.  Suppose that $\{\bar \a_1,\ldots,\bar \a_n\}\sub T$ and $\{r_1,\ldots,r_n\}\sub\bbbf$ with  $\sum_{i=1}^nr_i\bar \a_i=0.$ Take  $\{a_j\mid j\in J\}$ to be  a basis for $\bbbq$-vector space $\bbbf.$ For each $1\leq i\leq n ,$ suppose $\{r^j_i\mid j\in J\}\sub\bbbq$ is such that $r_i=\sum_{j\in J}r^j_ia_j.$ Then for each $\bar\a\in T,$ we have
\begin{eqnarray*}
0=\sum_{i=1}^nr_i\frac{2(\bar \a_i,\bar\a)^{\bar{}}}{(\bar\a,\bar\a)^{\bar{}}}=\sum_{i=1}^n\sum_{j\in J}r^j_ia_j\frac{2(\bar \a_i,\bar\a)^{\bar{}}}{(\bar\a,\bar\a)^{\bar{}}}=\sum_{j\in J}\sum_{i=1}^nr^j_i\frac{2(\bar \a_i,\bar\a)^{\bar{}}}{(\bar\a,\bar\a)^{\bar{}}}a_j.
\end{eqnarray*}
Since $\frac{2(\bar \a_i,\bar\a)^{\bar{}}}{(\bar\a,\bar\a)^{\bar{}}}\in\bbbz,$ we get that for each $j\in J$ and $\bar\a\in T,$ $$(\sum_{i=1}^nr^j_i\bar \a_i,\bar\a)^{\bar{}}=\sum_{i=1}^nr^j_i(\bar \a_i,\bar\a)^{\bar{}}=0.$$ So by (\ref{non}), $\sum_{i=1}^nr^j_i\bar \a_i=0$  for all $j\in J.$ But $T$ is $\bbbq$-linearly independent  and so $r_i^j=0$ for all $1\leq i\leq n$ and $ j\in J.$
This  means that  \begin{equation}\label{lin}
\parbox{2in}{$T$ is  $\bbbf$-linearly independent.}
 \end{equation}
Next suppose that $\bar R$ is of imaginary type and fix $\a^*\in\rim^\times.$ Using a modified version of the above argument together with \cite[Lem 3.14]{you2} (see also \cite[Lem. 3.21]{you2}), we get that \begin{equation}\label{lin2}
\parbox{2.5in}{$T\cup\{\bar{\a^*}\}$ is  $\bbbf$-linearly independent.}
\end{equation}
For each element $\a\in T,$ we fix a preimage $\dot\a\in R$ of $\a$ under $\bar{\;}$ and set
$$K:=\left\{\begin{array}{ll}
\{\dot\a\mid \a\in T\}&\hbox{if $\bar R$ is of real type,}\\
\{\dot\a\mid \a\in T\}\cup\{\a^*\}&\hbox{if $\bar R$ is of imaginary  type.}
\end{array}\right.$$
We have  using \cite[Pro. 3.14]{you2} together with (\ref{ref}) that $\bar \v=\hbox{span}_\bbbf \bar K.$ Therefore setting  $\dot \v:=\hbox{span}_\bbbf K$ and using (\ref{lin}) and (\ref{lin2}), we get that  $\v=\dot\v\op \v^0.$ We set $\dot R:=\{\dot\a\in\dot \v\mid \exists\sg\in\v^0, \dot \a+\sg\in R\},$ then  $\dot R$ is a locally finite root supersystem in its $\bbbz$-span isomorphic to $\bar R.$ Also since $K\sub R\cap\dot R,$ $-K\sub R\cap\dot R.$ So the subgroup $\w_K$ of the Weyl group of $R$  generated by the reflections based on  real roots of $K,$ we have
$$\w_K(\pm K)\sub R\cap \dot R\andd \pm\w_KK=\left\{\begin{array}{ll}(\dot R_{re})_{red}^\times&\hbox{if $\bar R$ is of real type,}\\
\dot R^\times& \hbox{if $\bar R$ is of imaginary  type.}\end{array}\right.$$%\begin{eqnarray*}
%\hbox{span}_\bbbq(1\ot \dot R)=\hbox{span}_\bbbq (1\ot S)=\bbbq\ot \hbox{span}_\bbbz S=\bbbq\ot \hbox{span}_\bbbz \dot R_{re}= \hbox{span}_\bbbq 1\ot \dot R_{re}.
%\end{eqnarray*}
We finally  set $S_{\dot \a}:=\{\sg\in \v^0\mid \dot\a+\sg\in R\}$ for $\dot\a\in \dot R.$ Then $R=\cup_{\dot\a\in\dot R}(\dot\a+S_{\dot\a})$ and  $$0\in S_{\dot\a} \hbox{ for }
\dot\a\in \left\{\begin{array}{ll}
(\dot R_{re})_{red} & \hbox{\small $\dot R$ is of real type,}\\
\dot R & \hbox{\small $\dot R$ is of imaginary  type.}
\end{array}\right.$$ Other assertions in the statement follow from  the same argument as in Claims 3,4 of the proof of Theorem 3.17 of \cite{you3}.\qed

\section{Root graded Lie superalgebras}

\begin{deft}
{\rm For a locally finite root supersystem  $R$ of type $X.$ Set
$$R_0:=\left\{\begin{array}{ll}
\{\a\in \rre\mid 2\a\not\in R\}\cup\{0\}&\hbox{if $X\neq BC(T,T')$}\\\\
\rre\setminus(\rre^2)_{sh}&\hbox{if $X= BC(T,T')$ and $\rre=\rre^1\op\rre^2$ }\end{array}\right.$$

and $$R_1:=R\setminus R_0.$$
We call elements of $R_0$ (resp. $R_1$) {\it even} (resp. {\it odd}) roots.}
\end{deft}

We note that for  a locally finite root supersystem $R,$ $R_0$ is a locally finite  root system.
\begin{deft}\label{root-graded-def}{\rm
 Suppose that $(A,\fm,R)$ is a locally finite root supersystem and $\Lam$ is an additive  abelian group. A Lie superalgebra $\LL=\LL_0\op\LL_1$ is called  an {\it $(R,\Lam)$-graded} Lie superalgebra if
\begin{itemize}
\item  the Lie superalgebra $\LL$ is equipped with  a $\la R\ra$-grading  $\LL=\op_{\a\in \la  R\ra}\LL^\a,$ that is
 \begin{itemize}
\item  $\LL_0$ as well as $\LL_1$ are $\la R\ra$-graded subspaces,
\item $[\LL^\a,\LL^\b]\sub\LL^{\a+\b}$ for all  $\a,\b\in\la R\ra,$
 \end{itemize}
 \item the support of $\LL$ with respect to the $\la R\ra$-grading is a subset of $R,$
\item $\LL^0=\sum_{\a\in  R\setminus\{0\}}[\LL^\a,\LL^{-\a}],$
\item the Lie superalgebra $\LL$ is equipped with a $\Lam$-grading $\LL=\op_{\lam\in\Lam}{}^\lam\LL$ which is compatible with the $\la R\ra$-grading on $\LL,$ that is
\begin{itemize}
\item  $\LL_0$ as well as $\LL_1$ are $\Lam$-graded subspaces,
\item $\LL^\a$ is a $\Lam$-graded subspace for each $\a\in R,$
\item $[{}^\lam\LL,{}^\mu\LL]\sub{}^{\lam+\mu}\LL$ for all  $\lam,\mu\in\Lam,$
\end{itemize}
\item there is a full subsystem $\Phi$ of $R$ such that for $0\neq\a\in \Phi,$ there are $0\neq e\in{}^0\LL\cap\LL^\a$ and $0\neq f\in {}^{0}\LL\cap\LL^{-\a}$ with $k_\a:=[e,f]\in \LL_0\setminus\{0\}$ and  for
$\b\in  R$ and $x\in\LL^\b,$ $[k_\a,x]=(\b,\a)x$ (we call $\{k_\a\mid \a\in \Phi\setminus\{0\}\}$ a {\it set of toral elements} and refer to $\Phi$ as  a {\it grading subsystem}).
\end{itemize}
}
\end{deft}
An $(R,\Lam)$-graded Lie superalgebra $\LL$ is called {\it fine} if for $i=0,1,$ the support $\LL_i$ with respect to the $\la R\ra$-grading is a subset of $R_i;$ also it is called {\it predivision} if for $\a\in R\setminus\{0\}$ and $\lam\in \Lam$ with ${}^\lam\LL^\a:={}^\lam\LL\cap\LL^{\a}\neq \{0\},$ there are $e\in{}^\lam\LL^\a$ and $f\in {}^{-\lam}\LL^{-\a}$ such that  $k:=[e,f]\in\LL_0\setminus\{0\}$ and for
$\b\in  R$ and $x\in\LL^\b,$ $[k,x]=(\b,\a)x.$ An $(R,\{0\})$-graded Lie superalgebra is called an {\it $R$-graded Lie superalgebra.}

\begin{lem}
Suppose that $(A,\fm,R)$ is a locally finite root supersystem and $\Lam$  an additive abelian group. If $\gg=\op_{\a\in R}\op_{\sg\in\Lam}{}^\sg\gg^\a$ is an $(R,\Lam)$-graded Lie superalgebra with a grading subsystem $\Phi,$ then so is  $\gg/Z(\gg).$ Moreover, if $\gg$ is predivision, then $\gg/Z(\gg)$ is also predivision.
\end{lem}
\pf Since $Z(\gg)$ inherits the gradings on $\gg,$ for $\a\in R$ and $\sg\in \Lam,$ we have $$\frac{{}^\sg\gg+Z(\gg)}{Z(\gg)}\cap\frac{\gg^\a+Z(\gg)}{Z(\gg)}=\frac{{}^\sg\gg^\a}{Z(\gg)}$$  and that $$\frac{\gg}{Z(\gg)}=\frac{\gg_{\bar 0}+Z(\gg)}{Z(\gg)}\op\frac{\gg_{\bar 1}+Z(\gg)}{Z(\gg)}=\bigoplus_{\a\in R,\sg\in \Lam}\frac{{}^\sg\gg^\a+Z(\gg)}{Z(\gg)}.$$ More precisely, $\gg/Z(\gg)$ is equipped with compatible $\la R\ra$ and $\Lam$-gradings. Now we prove that  $Z(\gg)\sub \gg^0.$ For this, we suppose $\a\in R\setminus\{0\}$ and show that $\gg^\a \cap Z(\gg)=\{0\}.$  If $\gg^\a=\{0\},$ there is nothing to prove, so suppose $\gg^\a\neq \{0\}.$ Since $\hbox{span}_\bbbq \Phi=\bbbq \ot_\bbbz R,$ for each $\b\in R,$ there is a nonzero integer $n$ with $n\b\in\hbox{span}_\bbbz \Phi.$ This together with the fact that the form $\fm$ is nondegenerate on $A=\hbox{span}_\bbbz R,$ guarantees the existence of an element $\gamma\in \Phi$ with $(\a,\gamma)\neq 0.$ Suppose that $k_\gamma$ to be a toral element of $\gg$ corresponding to $\gamma.$ For each  $0\neq x\in \gg^\a,$ we have $[k_\gamma,x]=(\a,\gamma)x\neq 0,$ so $x\not\in Z(\gg).$ This shows that $\gg^\a\cap Z(\gg)=\{0\}.$ To complete the proof, it is enough to show if
 $e\in{}^\lam\gg^\a$ and  $f\in {}^{-\lam}\gg^{-\a}$ for some $\a\in R\setminus\{0\}$ and $\lam\in \Lam$  with  $k:=[e,f]\in\gg_0\setminus\{0\}$ such that  $[k,x]=(\b,\a)x$ for
$\b\in  R, x\in\gg^\b,$  then $k\not\in Z(\gg).$  So consider $\a,\lam,e,f$ and $k$ as above. Since $\a\neq0,$ as before, there is $\b\in \Phi$ with $(\a,\b)\neq 0.$ Now for  $0\neq y\in{}^0\gg^\b,$ we have  $[k,y]=(\b,\a)y\neq 0.$ This shows that $k\not\in Z(\gg)$ and so we are done.
\qed
\begin{lem}
Suppose that $(\LL,\hh,\fm)$ is an extended affine Lie superalgebra with irreducible  root system $R.$ Keep the same notations as in Lemma  \ref{lem} and set $\Lam:=\la \cup_{\dot\a\in \dot R}S_{\dot\a}\ra,$ then  the  core $\LL_c$ of $\LL$ is a predivision $(\dot  R,\Lam)$-graded  Lie superalgebra.  Moreover, if   $R^0\sub R_0,$ then   for $i=0,1,$ the support $(\LL_c)_{\bar i}$ with respect to the $\la \dot R\ra$-grading is $\dot R_i.$
\end{lem}
\pf
We note that   for each root $\a\in R,$ $\LL^\a$ is a $\bbbz_2$-graded subspace, so $\LL_c$ is a $\bbbz_2$-graded subalgebra of $\LL.$ Moreover,  we have $$\LL_c=\sum_{\dot\a\in\dot R^\times,\sg\in S_{\dot\a}}\LL^{\dot\a+\sg}+\sum_{\
\dot\a\in\dot R^\times,\sg\in S_{\dot\a},\tau\in S_{-\dot\a}}[\LL^{\dot\a+\sg},\LL^{-\dot\a+\tau}].$$ Therefore, we have  $$\LL_c=\sum_{\dot\a\in \dot R}(\LL_c)^{\dot\a}=(\LL_c)_0\op(\LL_c)_1=\sum_{\sg\in\Lam}{}^\lam(\LL_c)$$ where
$$\begin{array}{l}
(\LL_c)^{\dot\a}=\sum_{\sg\in S_{\dot\a}}\LL^{\dot\a+\sg}\;\;\;\;\; (\dot\a\in\dot R^\times),\\
(\LL_c)^0=\sum_{\dot\a\in\dot R^\times}\sum_{\sg\in S_{\dot\a}}\sum_{\tau\in S_{-\dot\a}}[\LL^{\dot\a+\sg},\LL^{-\dot\a+\tau}],\\
(\LL_c)_{\bar 0}=\LL_{\bar 0}\cap\LL_c\andd(\LL_c)_{\bar 1}=\LL_{\bar 1}\cap\LL_c,\\
{}^\lam(\LL_c)=\sum_{\dot\a\in\dot R^\times}\LL^{\dot\a+\lam}+\sum_{\dot\a\in\dot R^\times}\sum_{\sg\in S_{\dot\a}}[\LL^{\dot\a+\sg},\LL^{-\dot\a+\lam-\sg}]\;\;\;\;\; (\lam\in \Lam).
\end{array}
$$
These define   compatible  $\la \dot R\ra$ and $\Lam$-gradings  on $\LL_c.$ Now set $$\dot \Phi:=\left\{\begin{array}{ll}
(\dot R_{re})_{red}& \hbox{if $\dot R$ is of real type}\\
\dot R& \hbox{if $\dot R$ is of imaginary type.}\end{array}\right.$$
We  know from Lemma \ref{lem} that  $\dot\Phi\sub R.$ Now  for $\dot\a\in \dot\Phi\setminus\{0\},$ since $\LL$ is  an extended affine Lie superalgebra, by \cite[Lem. 2.4]{you3}, there are $e\in \LL^{\dot\a}={}^0(\LL_c)^{\dot\a}$ and $f\in \LL^{-\dot\a}={}^0(\LL_c)^{-\dot\a}$ such that $[e,f]=t_{\dot\a}$ (we recall $t_{\dot\a}$ from Subsection \ref{extended}). Therefore,  for $x\in {}^\lam(\LL_c)^{\dot\b}\sub\LL^{\dot\b+\lam}$ ($\dot\b\in \dot R,\; \lam\in\Lam$), we have $$[t_{\dot\a},x]=(\dot\b+\lam)(t_{\dot\a})x=(t_{\dot\b+\lam},t_{\dot\a})x=(\dot\b+\lam,\dot\a)x=(\dot\b,\dot\a)x.$$  Now assume $R^0\sub R_0,$ then
using the same argument as in \cite[Pro. 2.14]{you3}, one gets that
\begin{equation}\label{com1}
\begin{array}{l}
\hbox{$\bullet$ if $\dot\a\in \dot R_{re}$ and $2\dot\a\not\in \dot R,$ then $\dot\a+S_{\dot\a}\sub R_0,$}\\
\hbox{$\bullet$ if $\dot\a\in \dot R_{re}^\times$ and $2\dot\a\in \dot R,$ then  $2\dot\a+S_{2\dot\a}\sub R_0,$}\\
\hbox{$\bullet$ if $\dot\a\in \dot R_{im}^\times,$ then $\dot\a+S_{\dot\a}\sub R_1,$} \\
\hbox{$\bullet$ if $\dot \a\in \dot R_{re}^\times$ and $\dot\a+\sg\in R_0$ for some  $\sg\in S_{\dot\a},$ then $\dot\a+\tau\not\in R_1$ for all $\tau\in S_{\dot\a}.$}
\end{array}
\end{equation} Now the result easily  follows if $\dot R$ is not of type $BC(T,T'), B(T,T'), B(1,T), B(T,1)$ and $G(1,2).$ So we just consider these mentioned types.
% We carry out the proof for these types as in the following cases. In what follows we use the natation as in Table ...
%
%Case 1. $\dot R$ is of type $BC(T,T'):$
% If $\d_p+\sg\in R_1$ for some $p\in T'$ and $\tau\in S_{\d_p},$ then by (\ref{com1}), $\d_p+\sg\in R_1$ for all $\sg\in S_{\d_p}.$ Let $\sg\in S_{\d_p}.$ Fix $i\in T,$  set $\a:=\d_p+\sg$ $e_{\a}\in\LL^{\a}_1$  and  pick $e_{-\a}\in\LL^{-\a}_1$ such that $[e_\a,e_{-\a}]=2t_\a/(\a,\a).$ Now suppose that $\zeta\in S_{\ep_i-\d_p}$ and consider $\gamma:=\ep_i-\d_p+\zeta.$ We note that $\gamma\in R_1$ and that  $[e_{-\a},\LL^{\gamma}]\sub\LL^{\ep_i-2\d_p-\sg+\zeta}=\{0\},$ so   we get
%\begin{eqnarray*}
%\{0\}\neq 2(\gamma,\a)/(\a,\a)\LL^\gamma&=&[2t_\a/(\a,\a),\LL^{\gamma}]\\
%&=&[[e_\a,e_{-\a}],\LL^{\gamma}]\\
%&=&[e_\a,[e_{-\a},\LL^{\gamma}]]+[e_{-\a},[e_{\a},\LL^{\gamma}]]=[e_\a,[e_{-\a},\LL^{\gamma}]].
%\end{eqnarray*}
%which implies that  $\{0\}\neq [e_{\a},\LL^{\gamma}]\sub\LL^{\e_i+\sg+\zeta}.$
%Therefore, $$\ep_i+\tau+\zeta\in R_0\;\; (i\in T,\sg\in S,\zeta\in F).$$ But  $F=-F$ and $S+F=S,$ so we have   $$\ep_i+\sg\in R_0\;\; (i\in T,\sg\in S).$$ Similarly, we get that $$\d_q+\eta\in R_1\;\; (q\in T',\eta\in S).$$
%A similar argument as above shows that if $\a=\d_p+\tau\in R_0,$ then
%$$\ep_i+\sg\in R_1\;\; (i\in T,\sg\in S_{1})\andd \d_q+\eta\in R_0\;\; (q\in T',\eta\in S_{2})$$ and so without loss of generality, we assume
%$$\cup_{\dot\a\in    ({{{\dot R}_{re}}^1})_{sh}}(\dot\a+S_{\dot\a})\sub R_0\andd\cup_{\dot\a\in    ({{{\dot R}_{re}}^2})_{sh}}(\dot\a+S_{\dot\a})\sub R_1.$$
From the classification table of Theorem \ref{classification II}, we  know that for types  $BC(T,T'), B(T,T'), B(1,T), B(T,1)$ and $G(1,2),$ $\dot R_{re}$ has two irreducible components $\dot R^1_{re}$ and $\dot R^2_{re}$ and that $\dot R_{im}^\times=(\dot R^1_{re})_{sh}+(\dot R^2_{re})_{sh}.$  We also recall from Lemma \ref{lem} that
  $S=S_{\dot\a},$ for all $\dot\a\in (\dot R_{re})_{sh},$  and $F=S_{\dot\b},$ for all $\dot\b\in \dot R_{lg}\cup\dot R_{im},$ satisfy   $S+F= S.$
 Considering (\ref{com1}), to complete the proof, we just need to show that if $\{i,j\}=\{1,2\},$  $\{r,s\}=\{0,1\}$ and $\dot\a+\sg\in R_r$  for some $\dot \a\in (\dot R^i_{re})_{sh}$ and $\sg\in S,$ then $$(\dot R^i_{re})_{sh}+S\sub R_r\andd (\dot R^j_{re})_{sh}+S\sub R_s.$$  So suppose that $\dot\a+\sg\in R_r$  for some $\dot \a\in (\dot R^i_{re})_{sh}$ and $\sg\in S.$   By (\ref{com1}), $\dot\a+\tau\in R_r$  for all $\tau\in S.$
Fix $\dot \b\in (\dot R^j_{re})_{sh}$ and $\tau\in S.$  Set $\a:=\dot\a+\tau\in R_r$ and pick $e_{\a}\in\LL^{\a}_r$  and $e_{-\a}\in\LL^{-\a}_r$ such that $[e_\a,e_{-\a}]=t_\a.$ We know that for each $\zeta\in F,$ $\gamma:=\dot\b-\dot\a+\zeta\in R_1$ and that  $[e_{-\a},\LL^{\gamma}]\sub\LL_{r+1}^{\dot\b-2\dot\a-\tau+\zeta}=\{0\},$ so   we get
\begin{eqnarray*}
\{0\}\neq (\gamma,\a)\LL^\gamma&=&[t_\a,\LL^{\gamma}]\\
&=&[[e_\a,e_{-\a}],\LL^{\gamma}]\\
&=&[e_\a,[e_{-\a},\LL^{\gamma}]]-(-1)^{r}[e_{-\a},[e_{\a},\LL^{\gamma}]]\\
&=&-(-1)^{r}[e_{-\a},[e_{\a},\LL^{\gamma}]]
\end{eqnarray*}
which implies that  $ \{0\}\neq[e_{\a},\LL^{\gamma}]\sub\LL_{r+1}^{\dot\b+\tau+\zeta}.$
Therefore, $$\dot\b+\tau+\zeta\in R_s\;\;\;\;\; (\tau\in S,\zeta\in F).$$ But  $S+F=S,$ so we have   $\dot\b+\eta\in R_s$ for $\eta\in S.$ This means that $$(\dot R^j_{re})_{sh}+S\sub R_s.$$ Finally using  the same  argument as above, we get that   $(\dot R^i_{re})_{sh}+S\sub R_r.$ This completes the proof.
%
%
%
%
% $\dot R$ is of type $B(T,T'),B(1,T), B(T,1), G(1,2):$
%It is enough to show that if $\dot\a\in \dot R_{re}^\times$ and $\sg\in S_{\dot\a}$ with $\dot\a+\sg,2(\dot\a+\sg)\in R,$ then  $\dot\a+\sg\in R_1.$ So suppose that $\dot\a\in \dot R_{re}^\times$ and $\sg\in S_{\dot\a}$ with $\dot\a+\sg,\in R.$ We know from the classification theorem that  there is  $\dot\b\in(\dot R_{re})_{sh}$ such that $\dot\b+S\sub R_0$ and $\pm\dot\a\pm\dot\b+F\sub \dot R_{im}\sub R_1.$ Suppose that $\sg\in S$ and set $\a:=\dot\b+\sg.$ Fix
%$e_{\a}\in\LL^{\a}_0$  and $e_{-\a}\in\LL^{-\a}_0$ such that $[e_\a,e_{-\a}]=2t_\a/(\a,\a).$ Now suppose that $\zeta\in F$ and consider $\gamma:=\dot\a-\dot\b+\zeta.$ We note that $\gamma\in R_1$ and that  $[e_{-\a},\LL^{\gamma}]\sub\LL^{\dot\a-2\dot\b-\sg+\zeta}=\{0\},$ so  if  $[e_\a,\LL^{\gamma}]=\{0\},$ we get
%\begin{eqnarray*}
%\{0\}\neq 2(\gamma,\a)/(\a,\a)\LL^\gamma&=&[2t_\a/(\a,\a),\LL^{\gamma}]\\
%&=&[[e_\a,e_{-\a}],\LL^{\gamma}]\\
%&=&[e_\a,[e_{-\a},\LL^{\gamma}]]+[e_{-\a},[e_{\a},\LL^{\gamma}]]=\{0\}
%\end{eqnarray*}
%which is a contradiction, so $\{0\}\neq [e_{\a},\LL^{\gamma}]\sub\LL^{\dot\a+\sg+\zeta}.$
%Therefore, $$\dot\a+\sg+\zeta\in R_1\;\; (i\in T,\tau\in S,\zeta\in F).$$ But  $F=-F$ and $S+F=S,$ so we have   $$\dot\a+\sg\in R_1\;\; (\sg\in S).$$
%
%
%
\qed

\begin{lem}\label{root-graded}
Suppose that $(\dot A,\fm^{}\dot{},\dot R)$ is a locally finite root supersystem, $\Lam$  a torsion free additive  abelian group and  $\gg=\gg_0\op\gg_1=\op_{\lam\in\Lam}\op_{\dot\a\in \dot R}{}^\lam\gg^{\dot \a}$  an  $(\dot R,\Lam)$-graded Lie superalgebra with a grading subsystem $\dot\Phi$ and a set  $\{k_{\dot\a}\mid \dot\a\in\dot\Phi\setminus\{0\}\}$  of toral elements.
 For $\dot\a\in\dot\Phi\setminus\{0\},$ fix $e_{\dot\a}\in{}^{0}\gg^{\dot\a}$ and $f_{\dot\a}\in{}^{0}\gg^{-\dot\a}$ such that $k_{\dot\a}=[e_{\dot\a},f_{\dot\a}]$ and take $T$ to be the linear span of  $\{k_{\dot\a}\mid \dot\a\in\dot\Phi\setminus\{0\}\}.$ Suppose that $\gg$ is equipped with an even nondegenerate invariant suppersymmetric  bilinear form $\fm.$

\smallskip

(i) Suppose that $\dot R=\op_{j\in J}\dot R^{(j)}$ is the decomposition of $\dot R$ into irreducible sub-supersystems. Suppose that $\{\dot\a_i\mid i\in I\}\sub \dot\Phi$ is such that $\{k_{\dot\a_i}\mid i\in I\}$ is a basis for $T.$ If $j\in J$ and  $\dot \gamma\in \dot\Phi^{(j)}:= \dot R^{(j)}\cap \dot\Phi,$ then  $k_{\dot\gamma}\in\hbox{span}_\bbbf \{k_{\dot\a_i}\mid i\in I,\dot\a_i\in \dot\Phi^{(j)}\}.$ Moreover, if  $\{r_{i}\mid i\in I\}\sub\bbbf$ with $k_{\dot\gamma}=\sum_{i\in I}r_{i}k_{\dot\a_i},$ we have $\dot\gamma=\sum_{i\in I}r_{i}{\dot\a_i}\in\bbbf\ot_{\bbbz}\dot A$ (here we identify $\dot A$ as a subset of $\bbbf\ot_\bbbz \dot A$); in particular,
for $\dot\b\in \dot R,$  $$\begin{array}{l}\tilde{\dot\b}:T \longrightarrow \bbbf\\
 k_{\dot\a}\mapsto(\dot\a,\dot\b\dot)\;\;\;(\dot\a\in \dot\Phi\setminus\{0\})
 \end{array}$$ is a well defined linear function.

\smallskip

(ii)  $\gg$ has a weight space decomposition with respect to $T$ with the   set of weights contained in $\{\tilde{\dot\b}\mid \dot\b\in \dot R\}.$
\smallskip

(iii) Suppose that $\gg$ is centerless. Assume that $\dot\gamma\in \dot R\setminus\{0\}$ and there are  $e\in \gg^{\dot\gamma}$ and  $f\in\gg^{-\dot\gamma}$ such that   $k:=[e,f]$ satisfies  $$[k,x]=(\dot\b,\dot\a\dot)x\;\; (x\in \gg^{\dot\b}).$$ If  $\{r,r_{\dot\a}\mid \dot\a\in \dot\Phi\setminus\{0\}\}\sub\bbbz\setminus\{0\}$ and  $r\dot\gamma=\sum_{\dot\a\in\dot\Phi\setminus\{0\}}r_{\dot\a}\dot\a, $ then $rk=\sum_{\dot\a\in\dot\Phi\setminus\{0\}}r_{\dot\a}k_{\dot\a};$ in particular, $k\in T$ and $(e,f)\neq 0.$

(iv) Suppose that $\gg$ is centerless. Assume that $\dot\gamma\in \dot R\setminus\{0\}$ and there are  $e,x\in \gg^{\dot\gamma}$ and  $f,y\in\gg^{-\dot\gamma}$ such that  $t:=[x,y]$ and $k:=[e,f]$ satisfy  $$[t,x]=(\dot\b,\dot\a\dot)x\andd [k,x]=(\dot\b,\dot\a\dot)x\;\; (x\in \gg^{\dot\b}).$$ Then $t=k$ and $(x,y)=(e,f).$
\end{lem}
\pf
$(i)$
The form  $\fm\dot{}$  induces the  $\bbbf$-bilinear form $\fm_\bbbf:(\bbbf\ot_\bbbz\dot A)\times(\bbbf\ot_\bbbz\dot  A)\longrightarrow \bbbf$ defined  by
$$(r\ot   a,s\ot   b)_\bbbf:=rs( a, b\dot);\;r,s\in\bbbf,\;a,b\in \dot A.$$ This is a nondegenerate symmetric bilinear form satisfying  $(\hbox{span}_\bbbz \dot R^{(i)},\hbox{span}_\bbbz \dot R^{(j)})_\bbbf=\{0\}$ for $i,j\in J$ with $i\neq j$ (see \cite[Lem. 3.21]{you2}).
Since   $\hbox{span}_\bbbq \dot\Phi=\bbbq \ot_\bbbz\dot R$ and  $\hbox{span}_\bbbz R=\op\hbox{span}_\bbbz\dot R^{(j)},$ we get that
$$
\hbox{span}_\bbbf  \dot\Phi^{(j)}=\hbox{span}_\bbbf \dot R^{(j)}.
$$
\begin{comment}
if $\a\in \hbox{span}_\bbbz \dot R^{(j)}\cap \op_{t\neq j}\hbox{span}_\bbbz \dot R^{(t)},$ then \begin{eqnarray*}
(\a,R)&=&(\a,\dot R^{(j)})\cup(\a,\cup_{t\neq j}\cup\dot R^{(t)})\\
&\sub&(\op_{t\neq j}\hbox{span}_\bbbz \dot R^{(t)},\dot R^{(j)})\cup (\hbox{span}_\bbbz \dot R^{(j)},\cup_{t\neq j}\dot R^{(t)})=0
\end{eqnarray*}
which contradicts the nondegeneracy of the form on $\dot A.$ Also if $\a\in \dot R^{(j)},$ then   $\a\in \hbox{span}_\bbbf 1\ot \dot\Phi,$ so there  is $x\in \hbox{span}_\bbbf 1\ot \dot\Phi^{(j)}, y\in  \sum_{j\neq t}\hbox{span}_\bbbf 1\ot \dot\Phi^{(t)}$ and $\a=x+y,$ so $\a-x=y,$ as above $(\a-x, R)_\bbbf=0,$
so $\a-x=0.$
\end{comment}
Suppose that $j\in J$ and $\dot\gamma\in \dot\Phi^{(j)}.$  Let $i_1,\ldots,i_n,$ $ j_1,\ldots,j_m\in I$ are such that  $\dot\a_{i_1},\ldots,\dot\a_{i_n}\in \dot \Phi^{(j)},$  $\dot\a_{j_1},\ldots,\dot\a_{j_m}\not\in  \dot \Phi^{(j)}$  and $k_{\dot\gamma}=r_1k_{\dot\a_{i_1}}+\cdots+r_nk_{\dot\a_{i_n}}+s_1k_{\dot\a_{j_1}}+\cdots+s_mk_{\dot\a_{j_m}}$ for some $r_1,\ldots,r_n,$ $s_1,\ldots,s_m\in\bbbf.$ For $\dot\b\in \dot\Phi^{(j)}\setminus\{0\},$ we have  $\gg^{\dot\b}\neq \{0\}$ and for $0\neq x\in \gg^{\dot\b},$ we have
\begin{eqnarray*}
(\dot\gamma,\dot\b)_\bbbf x&=&
(\dot\gamma,\dot\b\dot)x\\&=&[k_{\dot\gamma},x]\\
&=&[r_1k_{\dot\a_{i_1}}+\cdots+r_nk_{\dot\a_{i_n}}+s_1k_{\dot\a_{j_1}}+\cdots+s_mk_{\dot\a_{j_m}},x]\\
&=&r_1[k_{\dot\a_{i_1}},x]+\cdots+r_n[k_{\dot\a_{i_n}},x]+s_1[k_{\dot\a_{j_1}},x]+\cdots+s_m[k_{\dot\a_{j_m}},x]\\
&=&r_1(\dot\a_{i_1},\dot\b\dot)x+\cdots+r_n(\dot\a_{i_n},\dot\b\dot)x+s_1(\dot\a_{j_1},\dot\b\dot)x+\cdots+s_m(\dot\a_{j_m},\dot\b\dot)x\\
&=&r_1(\dot\a_{i_1},\dot\b\dot)x+\cdots+r_n(\dot\a_{i_n},\dot\b\dot)x\\
&=&(r_1\dot\a_{i_1}+\cdots+r_n\dot\a_{i_n},\dot\b)_\bbbf x
\end{eqnarray*}
This implies that $\dot\gamma=r_1\dot\a_{i_1}+\cdots+r_n\dot\a_{i_n}$ as the form $\fm_\bbbf$ on $\hbox{span}_\bbbf\dot R^{(j)}=\hbox{span}_\bbbf\dot\Phi^{(j)}$ is nondegenerate. Now we claim that $k_{\dot\gamma}=r_1k_{\dot\a_{i_1}}+\cdots+r_nk_{\dot\a_{i_n}}.$ To show this, it is enough to prove that $(k_{\dot\gamma}-(r_1k_{\dot\a_{i_1}}+\cdots+r_nk_{\dot\a_{i_n}}),k_{\dot\b})=0$ for all $\dot \b\in \dot\Phi\setminus\{0\}$ as the form is nondegenerate on $T.$ Assume $\dot\b\in \dot\Phi\setminus\{0\},$ then
\begin{eqnarray*}
(k_{\dot\gamma}-(r_1k_{\dot\a_{i_1}}+\cdots+r_nk_{\dot\a_{i_n}}),k_{\dot\b})&=&
(k_{\dot\gamma}-(r_1k_{\dot\a_{i_1}}+\cdots+r_nk_{\dot\a_{i_n}}),[e_{\dot\b},f_{\dot\b}])\\
&=&
([k_{\dot\gamma}-(r_1k_{\dot\a_{i_1}}+\cdots+r_nk_{\dot\a_{i_n}}),e_{\dot\b}],f_{\dot\b})\\
&=&(\dot\gamma-r_1\dot\a_{i_1}+\cdots+r_n\dot\a_{i_n},\dot\b)_\bbbf(e_{\dot\b},f_{\dot\b})\\
&=&0(e_{\dot\b},f_{\dot\b})=0.
\end{eqnarray*}  This completes the proof.
\begin{comment}
This shows that $k_{\dot\gamma}=r_1k_{\dot\a_{i_1}}+\cdots+r_nk_{\dot\a_{i_n}}.$
To complete the proof, we show that $\tilde{\dot\b}$ is a well-define linear map. For this, it is enough to consider $i_1,\ldots ,i_m\in I,$ $r_1,\ldots r_m\in \bbbf$ and  show if $\dot\gamma\in \dot R$ with $k_{\dot\gamma}=r_1k_{\dot\a_1}+\cdots + r_mk_{\dot\a_m}$
$\tilde{\dot\b}(k_{\dot\gamma})=r_1\tilde{\dot\b}(k_{\dot\a_1})+\cdots + r_m\tilde{\dot\b}(k_{\dot\a_m}).$  But we know  $\dot \gamma:=r_1{\dot\a_1}+\cdots + r_m{\dot\a_m},$ so we have
\begin{eqnarray*}\tilde{\dot\b}(k_{\dot \gamma})=(\dot\gamma,\dot\b\dot)=(\dot\gamma,\dot\b)_\bbbf&=&(r_1{\dot\a_1}+\cdots + r_m{\dot\a_m},\dot\b)_\bbbf\\
&=&r_1({\dot\a_1},\dot\b)_\bbbf+\cdots + r_m({\dot\a_m},\dot\b)_\bbbf\\
&=&r_1({\dot\a_1},\dot\b\dot)+\cdots + r_m({\dot\a_m},\dot\b\dot)\\
&=&r_1\tilde{\dot\b}(k_{\dot\a_1})+\cdots + r_m\tilde{\dot\b}(k_{\dot\a_m}).
\end{eqnarray*}
Also we note that if $\dot\eta,\dot\zeta\in \dot R$ with $k_{\dot\eta}=k_{\dot\zeta}=\sum_{i\in I} r_ik_{\dot\a_i},$ then $\dot\eta=\dot \zeta=\sum_{i\in I} r_i{\dot\a_i}\in\bbbf\otimes \dot A$ and so $\dot\eta=\dot\zeta$ so $(\dot\eta,\dot\b\dot)=(\dot\zeta,\dot\b\dot).$
\end{comment}

$(ii)$ % Suppose that $\dot\a,\dot\b\in\dot\Phi\setminus\{0\},$ then we have
%\begin{eqnarray*}
%[k_{\dot\a}, k_{\dot\b}]=[[e_{\dot\a},f_{\dot\a}],k_{\dot\b}]
%&=&[e_{\dot\a},[f_{\dot\b},k_{\dot\b}]]-(-1)^{|e_{\dot\a}|}[f_{\dot\b},[e_{\dot\a},k_{\dot\b}]]\\
%&=&(\dot\a,\dot\b\dot)[e_{\dot\a},f_{\dot\b}]+(-1)^{|e_{\dot\a}|}(\dot\a,\dot\b\dot)[f_{\dot\b},e_{\dot\a}]\\
%&=&0.
%\end{eqnarray*}
We know that $\gg=\op_{\dot\b\in \dot R}\gg^{\dot\b}.$ If $x\in \gg^{\dot\b}$ $(\dot\b\in\dot R),$ we have $[k_{\dot\a},x]=(\dot\b,\dot\a\dot)x=\tilde{\dot\b}(k_{\dot\a})x$ for all $\dot\a\in\dot\Phi\setminus\{0\}$ and so we have   $[t,x]=\tilde{\dot\b}(t)x$ for all $t\in T.$
So $\gg$ has a weight space decomposition $\gg=\op_{\dot\b\in\dot R}\gg^{(\tilde{\dot\b})}$ with respect to $T$ in which for $\dot\b\in\dot R,$ $\gg^{(\tilde{\dot\b})}=\gg^{\dot\b}.$

$(iii)$ We know $\gg=\sum_{\dot\b\in \dot R}\gg^{\dot\b}$ and that for all $a\in \gg^{\dot\b}$ ($\dot\b\in \dot R$),
\begin{eqnarray*}
[rk-\sum_{\dot\a\in \dot\Phi\setminus\{0\}}r_{\dot\a}k_{\dot\a},a]&=&
r(\dot\b,\dot\gamma\dot)a-\sum_{\dot\a\in\dot\Phi\setminus\{0\}}r_{\dot\a}(\dot\b,{\dot\a}\dot)a\\
&=&(\dot\b,r\dot\gamma-\sum_{\dot\a\in\dot\Phi\setminus\{0\}}r_{\dot\a}{\dot\a}\dot)a
=0.
\end{eqnarray*}  This means that  $rk-\sum_{\dot\a\in\dot\Phi\setminus\{0\}}r_{\dot\a}k_{\dot\a}$ is an element of the center of $\gg$ and so it is zero, i.e. $rk=\sum_{\dot\a\in\dot\Phi\setminus\{0\}}r_{\dot\a}k_{\dot\a};$ in particular, $k\in T.$ Now to the contrary, assume $(e,f)=0,$ then for each $\dot\a\in \dot\Phi\setminus\{0\},$ $$(k,k_{\dot\a})=([e,f],k_{\dot\a})=(e,[f,k_{\dot\a}])=(\dot\a,\dot\gamma\dot)(e,f)=0.$$ This contradicts the fact that the form on $T$ is nondegenerate.

$(iv)$ As $\gg$ is centerless, it is immediate that  $t=k.$ We shall show that  $(e,f)=(x,y).$ Since $\dot\gamma\neq 0,$ there is $\dot\a\in \dot\Phi$ with $(\dot\a,\dot\gamma\dot)\neq 0.$ Now we have
\begin{eqnarray*}
(e,f)(\dot\gamma,\dot\a\dot)=((\dot\gamma,\dot\a\dot)e,f)=([k_{\dot\a},e],f)=(k_{\dot\a},[e,f])=(k_{\dot\a},k)
&=&(k,k_{\dot\a})\\
&=&(k,[e_{\dot\a},f_{\dot\a}])\\
&=&([k,e_{\dot\a}],f_{\dot\a})\\
&=&(\dot\gamma,\dot\a\dot)(e_{\dot\a},f_{\dot\a}).
\end{eqnarray*}
This implies that $(e,f)=(e_{\dot\a},f_{\dot\a}).$ Similarly $(x,y)=(e_{\dot\a},f_{\dot\a}).$ This completes the proof.
 \qed

\begin{Thm}
Suppose that $(\dot A,\fm^{\dot{}},\dot R)$ is a locally finite root supersystem and $\Lam$ is a torsion free additive  abelian group. Suppose that $\gg=\op_{\lam\in\Lam}\op_{\dot\a\in \dot R}{}^\lam\gg^{\dot \a}$ is a centerless  $(\dot R,\Lam)$-graded Lie superalgebra, with a grading subsystem $\dot\Phi,$  equipped with an invariant  nondegenerate   even supersymmetric  bilinear form $\fm.$ Suppose that
 \begin{itemize}
\item for $\lam,\mu\in\Lam$ with $\lam+\mu\neq0,$ $({}^{\lam}\gg,{}^{\mu}\gg)=\{0\},$
\item the form is nondegenerate on the span of a set of toral elements of $\gg,$
\item for $\lam\in \Lam$ with ${}^\lam\gg_{\bar i}^{0}:=\gg_{\bar i}\cap{}^\lam\gg\cap\gg^{0}\neq \{0\}$ $(i=0,1),$ there are $e\in{}^\lam\gg_{\bar i}^{0}$ and $f\in {}^{-\lam}\gg_{\bar i}^{0}$ such that  $[e,f]_{_\gg}=0$ and $(e,f)\neq 0,$\\
\item for ${\dot\a}\in \dot R\setminus\{0\}$ and $\lam\in \Lam$ with ${}^\lam\gg_{\bar i}^{\dot\a}:=\gg_{\bar i}\cap{}^\lam\gg\cap\gg^{{\dot\a}}\neq \{0\}$ ($i=0,1$), there are $e\in{}^\lam\gg_{\bar i}^{\dot\a}$ and $f\in {}^{-\lam}\gg_{\bar i}^{-{\dot\a}}$ such that  $k:=[e,f]\in\gg_0\setminus\{0\}$ and for
$\dot\b\in \dot R$ and $x\in\gg^{\dot \b},$ $[k,x]=(\dot\b,{\dot\a}\dot)x,$
 \end{itemize} then $\gg$ is isomorphic to the  core of an extended affine Lie superalgebra modulo the center.
\end{Thm}
\pf Set $\v:=\bbbf\ot_\bbbz\Lam.$ Identify $\Lam$ with a subset of $\v$ and  fix a basis $\{\lam_i\mid i\in I\}\sub \Lam$ of $\v.$ Suppose that $\v^\dag$ is the restricted dual of $\v$ with respect to this basis. We suppose $\{d_i\mid i\in I\}$ is the corresponding basis for $\v^\dag.$
Consider $d_i$  ($i\in I$) as a derivation of $\gg$ mapping $x\in{}^\lam\gg$ to $d_i(\lam)x$ for all $\lam\in\Lam.$ Set  $$\LL:=\gg\op\v\op\v^\dag$$ and define

\begin{equation}\label{bracketloop}\begin{array}{l}
deg(v)=0;\;\; v\in\v\op\v^\dag\\
 \; -[x,d]=[d,x]:=d(\lambda) x,  \quad\quad (d\in\v^\dag,x\in \gg),\\
 \;[\LL,\v]=[\v,\LL]=[\v^\dag,\v^\dag]:=\{0\},\\
\;[x,y]:=[x,y]_{_\gg}+\sum_{i\in I}(d_i(x),y)\lambda_i, \quad\quad
(x,y\in \gg), \end{array}\end{equation} where by
$[\cdot,\cdot]_{\gg},$  we mean the Lie
bracket on ${\gg}.$
 We next extend the form on
$\gg$
  to a bilinear form on $\LL$ by
\begin{equation}\label{formloop}
\begin{array}{l}
(\v,  \v)=(\v^\dag,\v^\dag)=(\v,\gg)=(\v^\dag,\gg):=\{0\},\\
(v,d)=(d,v):=d(v),\quad  d\in \v^\dag, v\in \v.
\end{array}
\end{equation}
Then $\LL=\LL_{\bar 0}\op\LL_{\bar 1},$ where $\LL_{\bar 0}:=\gg_{\bar 0}\op\v\op\v^\dag$ and $\LL_{\bar 1}:=\gg_{\bar1},$ together with the Lie bracket $[\cdot,\cdot]$ is a Lie superalgebra and $\fm$ is an invariant  nondegenerate   even supersymmetric bilinear form.

For $\dot\a\in\dot\Phi\setminus\{0\},$  we fix $e_{\dot\a}\in{}^0\gg^{\dot\a}$ and $f_{\dot\a}\in{}^{0}\gg^{-\dot\a}$ such that  $k_{\dot\a}=[e_{\dot\a},f_{\dot\a}]$ and  that the form on $$T:=\hbox{span}_\bbbf\{k_{\dot\a}\mid \dot\a\in\dot\Phi\setminus\{0\}\}$$ is nondegenerate.  We next set $$\hh:=T\op\v\op\v^\dag$$ and note that the form restricted to $\hh$ is nondegenerate.  We identify $\hh^*$ with $T^*\op\v^*\op(\v^\dag)^*$ in the usual manner. We also consider $\lam\in\v$ as an element of $\hh^*$ by  $\lam(t+v+d)=d(\lam).$
We  know that $$\LL_{\bar 0}=\sum_{\lam\in\Lam}\sum_{\dot\a\in\dot R}{}^\lam\gg^{\dot\a}_{\bar 0}\op\v\op\v^\dag\andd \LL_{\bar 1}=\sum_{\lam\in\Lam}\sum_{\dot\a\in\dot R}{}^\lam\gg^{\dot\a}_{\bar 1}.$$
For $i\in\{0,1\},$ $t\in T,$ $v\in\v,$ $d\in\v^\dag,$ $\dot\b\in\dot R,$  $\lam\in\Lam$ and  $x\in {}^\lam\gg_{\bar i}^{\dot\b},$    we have using Lemma \ref{root-graded}($ii$) that
\begin{eqnarray*}&&[t+v+d, x]=[t,x]_{_\gg}+[d,x]=(\tilde{\dot\b}+\lam)(t+v+d)x,\\
&&[t+v+d, \v\op\v^\dag]=\{0\},
\end{eqnarray*}
so for $i=0,1,$ $\LL_{\bar i}$ has a  weight space decomposition with respect to $\hh$ with the set of weights $\{\tilde{\dot\b}+\lam\mid \dot\b\in \dot R,\lam\in\Lam,{}^\lam\gg_{\bar i}^{\dot\b}\neq\{0\}\}.$
Now suppose $i\in\{0,1\},$ $\dot\b\in \dot R,$ $\lam\in\Lam$ with $\tilde{\dot\b}+\lam\neq 0$ and $\LL_{\bar i}^{\tilde{\dot\b}+\lam}\neq\{0\}.$ So if $\dot\b\neq 0,$ there are $e\in{}^{\lam}\gg_{\bar i}^{\dot\b}$ and $f\in{}^{-\lam}\gg_{\bar i}^{-\dot\b}$ such that  $k:=[e,f]_{_\gg}\in\gg_0\setminus\{0\}$ and for $x\in \gg^{\dot\gamma},$ $[k,x]_{_\gg}=(\dot\b,\dot\gamma\dot)x.$
But since $\hbox{span}_\bbbq \dot \Phi=\bbbq\ot_\bbbz\hbox{span}_\bbbz \dot R,$  there is a nonzero integer $r\in\bbbz$ such that $r\dot\b=\sum_{\dot\a\in \dot \Phi\setminus\{0\}}r_{\dot\a}\dot\a.$ So $k=\frac{1}{r}\sum_{\dot\a\in\dot \Phi\setminus\{0\}}r_{\dot\a}k_{\dot\a}\in T$ by Lemma \ref{root-graded}.  This implies that  $[e,f]\in \hh\setminus\{0\}.$  Also if $\dot\b=0,$  take  $e\in{}^\lam\gg_{\bar i}^{0}$ and $f\in {}^{-\lam}\gg_{\bar i}^{0}$ such that  $[e,f]_{_\gg}=0$ and $(e,f)\neq 0,$ then $[e,f]\in\hh\setminus\{0\}.$ Therefore there are $e\in \LL_{\bar i}^{\tilde{\dot\b}+\lam}={}^\lam\gg_{\bar i}^{\dot\b}$ and $f\in \LL_{\bar i}^{-\tilde{\dot\b}-\lam}={}^{-\lam}\gg_{\bar i}^{-\dot\b}$ with $0\neq[e,f]\in \hh.$

Take $R$ to be the root system of $\LL$ with respect to $\hh$ and  suppose  $\dot\a\in \dot R$ and $\lam\in \Lam$ with $\tilde{\dot\a}+\lam\in R.$ If $\dot\a=0,$ then it is easy to see that $t_{\tilde{\dot\a}+\lam}=\lam$ in which $t_{\tilde{\dot\a}+\lam}$ as before is the unique element of  $\hh$ representing $\tilde{\dot\a}+\lam$ through the form $\fm.$ Also if $\dot\a\neq0,$ we fix $e\in{}^\lam\gg^{\dot\a}$ and $f\in{}^{-\lam}\gg^{-\dot\a}$ such that for $k:=[e,f]\in\gg_0\setminus\{0\}$ and for all $x\in\dot\gg^{\dot\b}$ ($\dot\b\in\dot R$) $[k,x]=(\dot\a,\dot\b\dot)x.$ Then considering Lemma \ref{root-graded}, it is easily verified that  $t_{\tilde{\dot\a}+\lam}=(e,f)^{-1}k+\lam.$ Now it follows that  $R^0=R\cap\Lam$ and $R^\times=\{\tilde{\dot\b}+\lam\mid \dot\b\in \dot R^\times,\lam\in\Lam,{}^\lam\gg^{\dot\b}\neq\{0\}\}.$
\begin{comment}
$R^0=\{\tilde{\dot\a}+\lam\mid {}^\lam\gg^{\dot\a}\neq \{0\};(k_{\dot\a},k_{\dot\b}})=\{0\}\forall \dot\b\in \dot \Phi\setminus\{0\}\}$
But if $\dot\a\neq 0$ and  $(k_{\dot\a},k_{\dot\b}})=\{0\}$ $\forall \dot\b\in \dot \Phi\setminus\{0\},$ we have $(e_{\dot\a},f_{\dot\a})(\dot\a,\dot\b\dot)=(e_{\dot\a},[f_{\dot\a},k_{\dot\b}}])=([e_{\dot\a},f_{\dot\a}],k_{\dot\b}})=\{0\},$ so $(\dot\a,\dot\b\dot)=0$ for all $\dot\b\in\dot R,$ we have $\dot\a=0$ which is a contradiction.
\end{comment}
We next show that $adx$ is locally nilpotent for  $x\in {}^\lam\gg^{\dot\a}=\LL^{\tilde{\dot\a}+\lam}$ ($\dot\a\in \dot R^\times$ and $\lam\in\Lam$ with  $\tilde{\dot\a}+\lam\in R$).  Let  $v\in \v$ and  $d\in\v^\dag,$ then  $adx(v)=0$ and if $\lam=0,$  $adx(d)=0.$ Next suppose that  $\lam\neq 0,$ then we have
\begin{eqnarray*}
(adx)^3(d)=-\lam(d)(adx)^2(x)&=&-\lam(d)[x,[x,x]]\\
&=&-\lam(d)[x,[x,x]_{_\gg}]\\&=&-\lam(d)([x,[x,x]_{_\gg}]_{_\gg}+\sum_{i\in I}(d_i(x),[x,x]_{_\gg})\lam_i\\
&=&-\lam(d)[x,[x,x]_{_\gg}]_{_\gg}\in\gg^{3\dot\a}=\{0\}.\end{eqnarray*}
Also for $\dot\b\in\dot R,$ $\mu\in\Lam$ and $y\in{}^\mu\gg^{\dot\b},$  since  $\dot R$ is a locally finite root supersystem, $\{k\dot\a+\dot\b\mid k\in\bbbz\}\cap\dot R$ is a finite set. Fix a positive integer $N$ such that for $m\geq N,$ $m\dot\a+\dot\b\not\in \dot R.$ If $\lam=0,$ we have $(adx)^N(y)=(ad_{_{\gg}}x)^N(y)\in \gg^{N\dot\a+\dot \b}=\{0\}$ in which $ad_{_\gg}$ denotes the adjoint representation of $\gg.$ If $\lam\neq0,$     we choose  a positive integer $n>N$ such that $n\lam+\mu\neq 0,$ then
\begin{eqnarray*}
(adx)^n(y)&=&(ad_{_\gg}x)^n(y)+\sum_{i\in I}(d_i(x),(ad_{_\gg}x)^{n-1}(y))\lam_i\\
&=&(ad_{_\gg}x)^n(y)\in \gg^{n\dot\a+\dot \b}=\{0\}.
\end{eqnarray*}
Therefore $adx$ is locally nilpotent. Thus  $(\LL,\fm,\hh)$ is an extended affine Lie superalgebra with root system $R=\{\tilde{\dot\b}+\lam\mid \dot\b\in \dot R,\lam\in\Lam,{}^\lam\gg^{\dot\b}\neq\{0\}\}.$ We now show that $\LL_c/Z(\LLc)$ is a Lie superalgebra  isomorphic to $\gg.$
We know that
$$\LLc=\sum_{\dot\a\in\dot R^\times,\lam\in\Lam}{}^{\lam}\gg^{\dot\a}+\sum_{\dot\a\in\dot R^\times}\sum_{\lam,\mu\in\Lam}[{}^{\lam}\gg^{\dot\a},{}^{\mu}\gg^{-\dot\a}]\sub\gg\op\v.$$
Take $\Pi$ to be the restriction of the canonical projection map $\LL\longrightarrow \gg$ to $\LLc$ with respect to the decomposition $\LL=\gg\op\v\op\v^\dag.$ Since $\gg^{\dot\a}\sub\LLc$ for all $\dot\a\in\dot R\setminus\{0\}$ and  $\gg^0=\sum_{\dot\a\in\dot R\setminus\{0\}}[\gg^{\dot\a},\dot\gg^{-\dot\a}]_{_\gg},$ $\Pi$ is surjective. Also if $x\in\gg$ and $v\in \v$ are such that $x+v\in Z(\LLc),$ then $[x+v,\gg^{\dot\a}]=\{0\}$ for all $\dot\a\in \dot R\setminus\{0\}.$ So $[x,\gg^{\dot\a}]_{_\gg}=\{0\}$ for all $\dot\a\in \dot R\setminus\{0\}.$  Then it follows that   $x\in Z(\gg)=\{0\}.$ Therefore $Z(\LLc)=\LLc\cap\v=ker\Pi.$ This implies that $\LLc$ is isomorphic to $\gg.$
 \qed

\begin{Example} {\rm Suppose that $\bbbf=\bbbc$ and take $I,J$ to be two disjoint  index sets with cardinal numbers greater than 2. We use the same notations as in Subsection \ref{bc}; in particular   $\mathfrak{u}$ is a vector superspace with a basis
$\{v_i\mid i\in I\cup\bar I\cup J\cup\bar J\cup\{0\}\}$ equipped with a supersymmetric bilinear form defined as in  (\ref{form1}). For $j,k\in \{0\}\cup I\cup \bar I\cup J\cup \bar J,$ consider $e_{j,k}$ and $\mathfrak{gl}$ as in (\ref{elementary2}) and (\ref{elemnatry2}) and for   $T=\sum_{j,k}r_{j,k}e_{j,k}\in \mathfrak{gl},$ set  $\bar T:=\sum_{j,k}\bar r_{j,k}e_{j,k}.$
%\begin{equation}\parbox{4in}{ for $i=0,1,$ if  $T\in \mathfrak{gl}_i,$ then $\bar T\in \mathfrak{gl}_i.$}\end{equation}
Now set $$\begin{array}{l}\LL_{\bar i}:=\{X\in \mathfrak{gl}_i\mid (Xv,w)=-(-1)^{|X||v|}(v,\bar Xw);\forall v,w\in \mathfrak{u}\};\;i=0,1\\
 \LL=\LL_{\bar 0}\op\LL_{\bar 1}\\
 \gg:=\LL\cap\hbox{span}_\bbbr\{e_{j,k}\mid, j,k\in\{0\}\cup I\cup \bar I\cup J\cup\bar J\},\\
 \hh:=\hbox{span}_\bbbr\{h_t:=e_{t,t}-e_{\bar t,\bar t},d_k:=e_{k, k}-e_{\bar k,\bar k}\mid  t\in I,\; k\in J\}.
 \end{array}$$
For $i\in I$ and $j\in J,$ define {\small$$\begin{array}{ll}
  \begin{array}{l}\ep_i:\hh\longrightarrow \bbbr\\
  h_t\mapsto \d_{i,t},\;\; d_k\mapsto 0,\\
\end{array}
  &
  \begin{array}{l}\d_j:\hh\longrightarrow \bbbr\\
  h_t\mapsto 0,\;\; d_k\mapsto \d_{j,k},
\end{array}
  \end{array}(t\in I,k\in J).$$}

One can see that with respect to $\hh,$ $\LL$ has a weight space decomposition $\LL=\op_{\a\in R}\LL^\a$ with the set of weights
{\small$$R:=\{\pm\ep_r,\pm(\ep_r\pm\ep_s),\pm\d_p,\pm(\d_p\pm\d_q),\pm(\ep_r\pm\d_p)\mid 1\leq r,s\leq m,\;1\leq p,q\leq n\}$$} and for $1\leq r\neq s\leq m$ and $1\leq p\neq q\leq n,$

{\footnotesize $$\begin{array}{ll}
\LL^{\ep_r}&=\hbox{span}_\bbbr\{i^{\a}(e_{r,0}-(-1)^\a e_{0,\bar r})\mid \a=0,1\},\vspace{2mm}\\
\LL^{-\ep_r}&=\hbox{span}_\bbbr\{i^{\a}(e_{\bar r,0}-(-1)^\a e_{0,r})\mid \a=0,1\},\vspace{2mm}\\
\LL^{2\ep_r}&=\hbox{span}_\bbbr  ie_{r,\bar r},\vspace{2mm}\\
\LL^{-2\ep_r}&=\hbox{span}_\bbbr ie_{\bar r,r},\vspace{2mm}\\
\LL^{2\d_p}&=\hbox{span}_\bbbr e_{p,\bar p},\vspace{2mm}\\
\LL^{-2\d_p}&=\hbox{span}_\bbbr e_{\bar p,p},\vspace{2mm}\\
\LL^{\ep_r+\ep_s}&=\hbox{span}_\bbbr\{i^{\a}(e_{r,\bar s}-(-1)^\a e_{s,\bar r})\mid \a=0,1\},\vspace{2mm}\\
\LL^{-\ep_r-\ep_s}&=\hbox{span}_\bbbr\{i^{\a}(e_{\bar r,s}-(-1)^\a e_{\bar s,r})\mid \a=0,1\},\vspace{2mm}\\
\LL^{\ep_r-\ep_s}&=\hbox{span}_\bbbr\{i^{\a}(e_{r,s}-(-1)^\a e_{\bar s,\bar r})\mid \a=0,1\},\vspace{2mm}\\
\LL^{\d_p+\d_q}&=\hbox{span}_\bbbr\{i^\a(e_{p,\bar q}-(-1)^{\a+1}e_{q,\bar p})\mid \a=0,1\},\vspace{2mm}\\
\LL^{-\d_p-\d_q}&=\hbox{span}_\bbbr\{i^\a(e_{\bar p,q}-(-1)^{\a+1}e_{\bar q,p})\mid \a=0,1\},\vspace{2mm}\\
\LL^{\d_p-\d_q}&=\hbox{span}_\bbbr\{i^\a(e_{p,q}+(-1)^{\a+1}e_{\bar q,\bar p})\mid \a=0,1\},\vspace{2mm}\\
\LL^{\d_p}&=\hbox{span}_\bbbr\{i^\a(e_{0,\bar p}+(-1)^{\a+1}e_{p,0})\mid \a=0,1\},\vspace{2mm}\\
\LL^{-\d_p}&=\hbox{span}_\bbbr\{i^\a(e_{0,p}-(-1)^{\a+1}e_{\bar p,0})\mid \a=0,1\},\vspace{2mm}\\
\LL^{\ep_r+\d_p}&=\hbox{span}_\bbbr\{i^\a(e_{r,\bar p}+(-1)^{\a+1}e_{p,\bar r})\mid \a=0,1\},\vspace{2mm}\\
\LL^{-\ep_r-\d_p}&=\hbox{span}_\bbbr\{i^\a(e_{\bar r,p}-(-1)^{\a+1}e_{\bar p,r})\mid \a=0,1\},\vspace{2mm}\\
\LL^{\ep_r-\d_p}&=\hbox{span}_\bbbr\{i^\a(e_{r,p}-(-1)^{\a+1}e_{\bar p,\bar r})\mid \a=0,1\},\vspace{2mm}\\
\LL^{-\ep_r+\d_p}&=\hbox{span}_\bbbr\{i^\a(e_{\bar r,\bar p}+(-1)^{\a+1}e_{p,r})\mid \a=0,1\}.
\end{array}$$}
It is easy to see that the Lie superalgebra $\fl:=\sum_{\a\in R^\times}\LL^\a+\sum_{\a\in R^\times}[\LL^\a,\LL^{-\a}]$ is a $B(I,J)$-graded Lie superalgebra with grading subsystem $(R_{re})_{red}$.}\end{Example}

\section{$BC(I,J)$-graded Lie superalgebras}
In this section we illustrate the structure of Lie superalgebras graded by the  locally finite  root supersystem $BC(I,J).$
Throughout this section, we use the same notations as in  Section \ref{bc}; in particular, we recall $\fg$ as well as $\mathfrak{s}$ from (\ref{gs}), $\D_\mathfrak{u}$ from (\ref{delta}) and $R,\D_\mathfrak{s}$ from (\ref{root}).
We set    $$\Psi:=\{\pm\ep_i\pm\ep_j,\pm\ep_i,\pm\d_p\pm\d_q,\pm\ep_i,\pm\d_p\mid i\in I,j\in J\}=R\cup\{\pm2\ep_i\mid i\in I\};$$ in fact $\Psi$  is a locally finite root supersystem of type $BC(I,J).$
Suppose that $\fl$ is a Lie superalgebra such that
\begin{equation}\label{root-gra}
\begin{array}{l}
\bullet
\hbox{ $\fl$ contains  $\fg$ as a subalgebra,}\\
\bullet \hbox{ $\fl$ is equipped with a weight space decomposition $\fl=\op_{\a\in \Psi}\fl^\a,$ with respect to $\fh,$} \\
\bullet \hbox{ $\fl^0=\sum_{\a\in \Psi\setminus\{0\}}[\fl^\a,\fl^{-\a}].$}
\end{array}
\end{equation}
It is easy to see that $\fl$ is a $\Psi$-graded Lie superalgebra with $R$ as its grading subsystem. One knows that (\ref{root-gra}) is just a generalization of the notion of root graded Lie superalgebra in the sense of \cite{BE1} by switching from  finite root supersystems  to locally finite root supersystems.
 In this section, we study the structure of a Lie superalgebra $\fl$ satisfying (\ref{root-gra}).
Throughout this section we suppose $\bbbf$ is an algebraically closed field of characteristic zero; we also make a convention that    for a map $f$ defined on a set $X,$ by $x^f,$ for $x\in X,$ we mean the image of $x$ under $f.$

\subsection{Some Conventions}\label{convention}
Suppose that $\fa$ is an associative superalgebra and $\eta$ is a {\it superinvolution}  of   $\fa,$ i.e.,  $\eta$ is  an even linear map with $\eta^2=id$ and  $\eta(ab)=(-1)^{|a||b|}\eta(b)\eta(a)$ for all $a,b\in\fa.$ Next we assume  $\cc$ is an associative $\fa$-supermodule and  $\chi:\cc\times\cc\longrightarrow \fa$ is  a  {\it superhermitian} $\fa$-form of $\cc,$ in the sense that $\chi$ is an even bilinear form satisfying $$\chi(x,y)^\eta=(-1)^{|x||y|}\chi(y,x)\andd \chi(ax,y)=a\chi(x,y)$$ for all $x,y\in \cc$ and $a\in\fa.$ Then $\fb:=\fb(\fa,\cc):=\fa\op\cc$ together with $$(\a+c)(\a'+c')=(\a\cdot\a'+2\chi(c,c'))+(\a\cdot c'+(-1)^{|\a'||c|}(\a')^\eta\cdot c)$$ is a superalgebra. We set $$\aa:=\{\a\in\fa\mid \a^\eta=\a\}\andd \bb:=\{\a\in\fa\mid \a^\eta=-\a\}$$ and note that $$\fa=\aa\op \bb.$$ We next define {\small $$\begin{array}{ll}
\diamond:\cc\times\cc\longrightarrow \aa &
(c, c')\mapsto\frac{1}{2}(\chi(c,c')+(-1)^{|c||c'|}\chi( c', c))\\\\
\heart:\cc\times\cc\longrightarrow \bb&
(c, c')\mapsto\frac{1}{2}(\chi(c,c')-(-1)^{|c||c'|}\chi(c', c)).
\end{array}
$$} Finally for  $\b_1=a_1+b_1+c_1,\b_2=a_2+b_2+c_3\in\aa\op\bb\op\cc,$  we set
   \begin{equation}\label{conv}\b^*_{\b_1,\b_2}:=[a_1,a_2]+[b_1,b_2]+2c_1\heart c_2,\; \b_1^*:=c_1\andd \b_2^*:=c_2.\end{equation}

\bigskip

\subsection{Structure Theorem}
We suppose $|I|,|J|>4$ and  fix a nonempty subset $I_0$ of $I$ of finite cardinal number $m>3$ as well as  a subset  $J_0$ of $ J$ of finite  cardinal number $n>3$ and set
\begin{equation}\label{m,n}\begin{array}{l}\Phi:=\{\pm\ep_i\pm\ep_j,\pm\ep_i,\pm\d_p\pm\d_q,\pm\d_p,\pm\ep_i\pm\d_p\mid i\in I_0,j\in J_0\},\\
 R^{^{m,n}}:=\Phi\setminus\{\pm2\ep_i\mid i\in I_0\}.\end{array}\end{equation}
  We next define  the linear endomorphism  $\hbox{id}_{m,n}$ on $\mathfrak{u}$ by
$$\begin{array}{c}\hbox{id}_{m,n}:\mathfrak{u}\longrightarrow \mathfrak{u}\\
v_0\mapsto v_0,\;v_i\mapsto v_i,\;v_{\bar i}\mapsto v_{\bar i},\;v_{j}\mapsto 0,\;v_{\bar j}\mapsto 0,\\
(\hbox{\small$i\in I_0\cup J_0,\; j\in I\cup J\setminus (I_0\cup J_0)$})
\end{array}$$
and  for $u,v\in\mathfrak{u}$ and $x,y\in\fg\cup\mathfrak{s},$ define
$$\begin{array}{l}\;[u, v]
:\mathfrak{u}\longrightarrow\mathfrak{u};\;w\mapsto(v,w)u+(-1)^{|u||v|}(u,w)v-\frac{2(u,v)}{2m+1-2n
}\hbox{id}_{m,n}(w);\;\;w\in\mathfrak{u},\vspace{2mm}\\
\;u\circ v:\mathfrak{u}\longrightarrow\mathfrak{u};\;w\mapsto(v,w)u-(-1)^{|u||v|}(u,w)v;\;\;w\in\mathfrak{u},\\
\;x\circ y:=xy+(-1)^{|x||y|}yx-\frac{2str(xy)}{2m+1-2n}\hbox{id}_{m,n}.
\end{array}$$

\begin{Thm}
\label{main}
Suppose that $\mathfrak{L}$ is a Lie superalgebra satisfying the following:
$$
\begin{array}{l}
\bullet
\hbox{ $\mathfrak{L}$ contains  $\fg$ as a subalgebra,}\\
\bullet \hbox{ $\mathfrak{L}$ is equipped with a weight space decomposition $\mathfrak{L}=\op_{\a\in \Psi}\mathfrak{L}^\a,$ with respect to $\fh,$} \\
\bullet \hbox{ $\mathfrak{L}^0=\sum_{\a\in \Psi\setminus\{0\}}[\mathfrak{L}^\a,\mathfrak{L}^{-\a}].$}
\end{array}
$$
 Then there are a subsuperalgebra $\dd$ of $\fl,$ superspaces $\aa,\bb,\cc,$ even bilinear maps  $$\begin{array}{llll}
\cdot:\fa\times\fa\longrightarrow \fa&
\cdot:\fa\times\cc\longrightarrow \cc,&
\chi:\cc\times\cc\longrightarrow \fa,&
\la\cdot,\cdot\ra:\fb\times\fb\longrightarrow \dd
\end{array}$$
in which $\fb:=\fa\op\cc,$ and linear maps $$\begin{array}{lc}\eta:\fa\longrightarrow \fa\andd&
\phi:\dd\longrightarrow End(\fb)
\end{array}$$
 such that
$(\fa,\cdot)$ is an associative superalgebra, $(\cc,\cdot)$ is an associative $\fa$-module, $\eta$ is a superinvolution and $\chi$ is a superhermitian $\fa$-form with the following properties:

\begin{itemize}
\item $\la\cdot,\cdot\ra$ is surjective, supersymmetric and  satisfy $\la \aa,\bb\ra=\la\aa,\cc\ra=\la\bb,\cc\ra=\{0\},$ where $\aa$ and $\bb$ are respectively fixed and skew-fixed points of $\aa$ under the map $\eta,$ \\
\item considering  the superalgebraic structure on  $\fb$ as constructed in Subsection \ref{convention},  for each $d\in\dd,$ we have $\phi(d)$ is a superderivation of $\fb;$ i.e., $\dd$ acts on $\fb$ as superderivations, \\
\item $d\aa\sub \aa,$ $d\bb\sub \bb$ and $d\cc\sub \cc,$ for all $d\in\dd,$\\
\item $[d,\la\b,\b'\ra]=\la d\b,\b'\ra+(-1)^{|d||b'|}\la \b,d\b'\ra,$\\
\item $\sum_\circlearrowleft(-1)^{|\b_1||\b_3|}\la\b_1,\b_2\b_3\ra=0,$\\
\item $\la\a,\a'\ra\a''=\frac{1}{2(2m+1-2n)}[[\a,\a']-[\a,\a']^\eta,\a''],$\\
\item $\la\a,\a'\ra c=\frac{1}{2(2m+1-2n)}([\a,\a']-[\a,\a']^\eta)c,$\\
\item $\la c,c'\ra \a=\frac{1}{2m+1-2n}[\chi(c,c')-\chi(c,c')^\eta,\a],$\\
\item $\la c,c'\ra c''=\frac{1}{2m+1-2n}(\chi(c,c')-\chi(c,c')^\eta)\cdot c''+(-1)^{|c|(|c'|+|c''|)}\chi(c',c'')^\eta\cdot c-(-1)^{|c'||c''|}\chi(c,c'')^\eta\cdot c'.$
\end{itemize}
 Moreover, we have the following:

(i) There are subsuperspaces $\mathfrak{L}^1,\mathfrak{L}^2$ and $\mathfrak{L}^3$ of $\mathfrak{L}$ isomorphic to $\fg\ot \aa,\mathfrak{s}\ot \bb $ and $\mathfrak{u}\ot\cc$ respectively such that  $$\mathfrak{L}=(\mathfrak{L}^1\op\mathfrak{L}^2\op\mathfrak{L}^3)+ \dd.$$ Furthermore, if either   $|I|=m$ and $ |J|=n$ or $I\cup J$ is an infinite set, we have
$$\mathfrak{L}=\mathfrak{L}^1\op\mathfrak{L}^2\op\mathfrak{L}^3\op \dd,$$ more precisely, in these cases $\mathfrak{L}$  can be identified with $$(\fg\ot \aa)\op(\mathfrak{s}\ot \bb)\op(\mathfrak{u}\ot\cc)\op\dd.$$

(ii) Identify $\mathfrak{L}^1,\mathfrak{L}^2$ and $\mathfrak{L}^3$ with $\fg\ot \aa,\mathfrak{s}\ot \bb $ and $\mathfrak{u}\ot\cc$ respectively, the Lie bracket on $\mathfrak{L}$ is given by the following:
{\small \begin{equation}\label{prepro}
\begin{array}{ll}
\;[x\ot a,y\ot a']&=(-1)^{|a||y|}([x,y]\ot\frac{1}{2}(a\circ a')+ (x\circ y)\ot\frac{1}{2}[a,a']+str(xy)\la a,a'\ra),\vspace{1mm}\\
\;[x\ot a,e\ot b]&= (-1)^{|a||e|}((x\circ e)\ot\frac{1}{2}[a,b]+[x,e]\ot\frac{1}{2}(a\circ b)),\vspace{1mm}\\
\;[e\ot b,f\ot b']&=(-1)^{|b||f|}([e,f]\ot\frac{1}{2}(b\circ b')+ (e\circ f)\ot\frac{1}{2}[b,b']+str(ef)\la b,b'\ra),\vspace{1mm}\\
\;[x\ot a,u\ot c]&=(-1)^{|a||u|}xu\ot a\cdot c,\vspace{1mm}\\
\;[e\ot b,u\ot c]&=(-1)^{|b||u|}eu\ot b\cdot c,\vspace{1mm}\\
\;[u\ot c,v\ot c']&=(-1)^{|c||v|}((u\circ v)\ot (c\diamond c')+ [u, v]\ot (c\heart c')+(u,v)\la c,c'\ra)\vspace{1mm}\\
~[\la\b_1,\b_2\ra,\la\b'_1,\b'_2\ra]&=\la\la\b_1,\b_2\ra\b'_1,\b'_2\ra+(-1)^{(|\b_1|+|\b_2|)|\b'_1|} \la\b'_1,\la\b_1,\b_2\ra\b'_2\ra,\vspace{1mm}\\
~[\la\b_1,\b_2\ra,x\ot a]&=\frac{(-1)^{|\b_1||x|+|\b_2||x|}}{2(2m+1-2n)}([\hbox{id}_{m,n},x]\ot (\b^*_{\b_1,\b_2}\circ a)+
(\hbox{id}_{m,n}\circ x)\ot [\b^*_{\b_1,\b_2}, a])\\\\
~[\la\b_1,\b_2\ra,e\ot b]&=\frac{(-1)^{|\b_1||x|+|\b_2||x|}}{2(2m+1-2n)}([\hbox{id}_{m,n},e]\ot (\b^*_{\b_1,\b_2}\circ b)+
(\hbox{id}_{m,n}\circ e)\ot [\b^*_{\b_1,\b_2}, b])\\
&-\frac{1}{2m+1-2n}str(\hbox{id}_{m,n}e)\la [b_1,b_2],b\ra\\\\
~[\la\b_1,\b_2\ra,u\ot c]&=\frac{(-1)^{|\b_1||u|+|\b_2||u|}}{2m+1-2n}(\hbox{id}_{m,n}u\ot \b^*_{\b_1,\b_2}c)
+(-1)^{|\b_1||u|+|\b_2||u|}u\ot \\&
((-1)^{|\b_1^*||\b_2^*|+|\b_1^*||c|}\chi(\b^*_2,c)^\eta \b^*_1-(-1)^{|\b^*_2||c|}\chi(\b^*_1,c)^\eta \b^*_2).
\end{array}
 \end{equation}
}

\end{Thm}

\begin{rem}
{\rm   We mention that if $|I|=m$ and $|J|=n,$ then the last three Lie brackets in the above display will be converted to the following ones:

  \begin{equation}
\begin{array}{ll}
~[\la\b_1,\b_2\ra,x\ot a]&=(-1)^{|\b_1||x|+|\b_2||x|}x\ot \la\b_1,\b_2\ra a\\\\
~[\la\b_1,\b_2\ra,e\ot b]&=(-1)^{|\b_1||e|+|\b_2||e|}e\ot \la\b_1,\b_2\ra b\\\\
~[\la\b_1,\b_2\ra,u\ot c]&=(-1)^{|\b_1||u|+|\b_2||u|}u\ot \la\b_1,\b_2\ra c.
\end{array}
 \end{equation}
\hfill$\diamondsuit$}
\end{rem}
To prove Theorem \ref{main}, we first  carry out the proof  for the case that  $|I|,|J|<\infty;$ at the first step, we suppose $I=I_0$ and $J=J_0.$
\subsubsection{Finite Case-The First Step}
In this subsection, we  assume  $I=I_0,$ $J=J_0$ and  that $\mathfrak{L}$ is a Lie superalgebra satisfying (\ref{root-gra}).
Consider $\mathfrak{L}$ as a $\fg$-module via the adjoint representation, then  $\mathfrak{L}$ is a locally finite $\fg$-module, i.e., any finite subset of $\mathfrak{L}$ generates a finite dimensional $\fg$-submodule (see \cite[Lem. 2.2]{BE3}).
Therefore, it is a summation of finite dimensional $\fg$-submodules. Using Corollary \ref{cor}, $\mathfrak{L}$ is completely reducible such that each of its irreducible components is either isomorphic to one of $\fg$-modules $\fg,$ $\frak{u},$ $\frak{s}$ or it is a trivial $\fg$-module.
Now collecting  the  isomorphic $\fg$-submodules  of the same parity,  we may assume that as a vector space, $\mathfrak{L}$ is isomorphic to $$(\fg\ot\aa_0)\op(\fg\ot\aa_1)\op(\mathfrak{s}\ot\bb_0)\op(\mathfrak{s}\ot\bb_1)\op(\frak{u}\ot\cc_0)\op(\frak{u}\ot\cc_1)\op\dd;$$
%We regard  $\fg\ot\aa_0$ as the summation of $\fg$-submodules isomorphic to $\fg$ of parity zero and $\fg\ot\aa_1$ as the summation of $\fg$-submodules isomorphic to $\fg$ of parity  $1$ and so on.
in which $\aa_0,$ $\aa_1,$ $\bb_0,$ $\bb_1,$ $\cc_0$ and $\cc_1$ are vector spaces and $\dd$ is the centralizer of  $\gg$ in $\mathfrak{L};$ in particular,  $\dd$ is a subsuperalgebra of $\LL.$
Setting $$\aa:=\aa_0\op\aa_1,\;\bb:=\bb_0\op\bb_1,\;\cc:=\cc_0\op\cc_1,\;\dd:=\dd_0\op\dd_1,$$ we can consider
$$\fl=(\fg\ot\aa)\op(\mathfrak{s}\ot\bb)\op(\frak{u}\ot\cc)\op\dd.$$
Now using the same argument as in \cite[$\S$ 5]{BE3}, one can see that the Lie superalgebraic structure on $\LL$ induces a superalgebraic structure on $\aa\op\bb\op\cc$ and that the stated properties  in  Theorem \ref{main} are fulfilled.

\subsubsection{Finite Case-Compatibility of Subsuperalgebras}\label{compatibility}
Throughout this subsection, we assume that $m'>m,$ $n'>n,$ $I=\{1,\ldots,m'\},J=\{1,\ldots,n'\}$  and $I_0=\{1,\ldots,m\},J_0=\{1,\ldots,n\}.$  We also assume   $$\fl:=\sum_{\a\in \Psi\setminus\{0\}}\fl^\a\op \sum_{\a\in \Psi\setminus\{0\}}[\fl^\a,\fl^{-\a}]$$ is a Lie superalgebra satisfying (\ref{root-gra}). Consider $\fl$ as a $\fg$-module. As in the previous subsection, it follows  from Corollary \ref{cor} that $\fl$ is decomposed into irreducible submodules, more precisely, \begin{equation}\label{weight1}\fl=\sum_{i\in I} \fg^{(i)}\op \sum_{j\in J} \v^{(j)}\op\sum_{t\in T} \mathfrak{s}^{(t)}\op E\end{equation} in which $\fg^{(i)}$ is isomorphic to $\fg,$ $\v^{(j)}$ is isomorphic to $\mathfrak{u},$ $\mathfrak{s}^{(t)}$ is isomorphic to $\mathfrak{s}$ for all $i\in I, j\in J,t\in T$ and $E$ is a trivial $\fg$-module. As before,  collecting  the  isomorphic $\fg$-submodules,  we   may assume
\begin{equation}\label{final}\fl=(\fg\ot\aa)\op(\mathfrak{s}\ot\bb)\op(\frak{u}\ot\cc)\op E,\end{equation} in which
$\aa,\bb,\cc$ are vector superspaces.
 We recall (\ref{m,n}) and use Proposition \ref{com} to get that  $$\gg:=\sum_{\a\in R^{^{m,n}}\setminus\{0\}} \fg^\a+\sum_{\a\in R^{^{m,n}}\setminus\{0\}}[\fg^\a,\fg^{-\a}]$$ is a subsuperalgebra of $\fg$ isomorphic to $\mathfrak{osp}(m,n)$ with Cartan subalgebra $$\hh=\sum_{\a\in R^{^{m,n}}\setminus\{0\}}[\fg^\a,\fg^{-\a}]$$ and root system $R^{^{m,n}}.$
Consider (\ref{m,n}) and set $$\LL:=\sum_{\a\in \Phi\setminus\{0\}}\fl^\a\op \sum_{\a\in \Phi\setminus\{0\}}[\fl^\a,\fl^{-\a}].$$ It is easy to see that   $\LL$ has a weight space decomposition $\LL=\sum_{\a\in \Phi}\LL^\a$ with respect to $\hh$ with $$\LL^\a:=\left\{\begin{array}{ll}\fl^\a&\a\in \Phi\setminus\{0\}\\
\sum_{\a\in \Phi\setminus\{0\}}[\fl^\a,\fl^{-\a}]&\a=0.\end{array}\right.$$ This in particular implies that \begin{equation}\label{weight}\LL^\a=\sum_{i\in I} (\fg^{(i)})^\a\op\sum_{j\in J}(\v^{(j)})^\a\op\sum_{t\in T}(\mathfrak{s}^{(t)})^\a\;\;\;(\a\in \Phi\setminus\{0\}).
\end{equation}
 Moreover, setting $$\D_{1}:=\hbox{span}_\bbbz R^{^{m,n}}\cap \D_\frak{u}\andd \D_2:=\hbox{span}_\bbbz R^{^{m,n}}\cap \D_{\frak{s}},$$ and using Proposition  \ref{com}, we have the following $\gg$-modules $$\begin{array}{ll}\hbox{ $\gg^{(i)}:=\sum_{\a\in R^{^{m,n}}\setminus\{0\}}(\fg^{(i)})^{\a}+\sum_{\a\in R^{^{m,n}}\setminus\{0\}}[\fg^\a,(\fg^{(i)})^{-\a}],$}\\\\
\hbox{ $\u^{(j)}:=\sum_{\a\in \D_1\setminus\{0\}}(\v^{(j)})^\a+\sum_{\a\in \D_1\setminus\{0\}}[\fg^\a,(\v^{(j)})^{-\a}],$ }\\\\
\hbox{ $\ss^{(t)}:=\sum_{\a\in  \D_2\setminus\{0\}}(\mathfrak{s}^{(t)})^{\a}+\sum_{\a\in \D_2\setminus\{0\}}[\fg^\a,(\mathfrak{s}^{(t)})^{-\a}]$}\end{array}$$ which are respectively  isomorphic to $\gg,$ to the natural module  $\u$ of $\gg$ and to  the second natural module $\ss$ of $\gg.$
Also it is immediate that $$
\begin{array}{l}
1)\hbox{ $\LL$ contains  $\gg$ as a subalgebra,}\\
2)\hbox{ $\LL$ is equipped with a weight space decomposition $\LL=$\hbox{\tiny$\displaystyle{\bigoplus_{\a\in \Phi}}$}$\LL^\a,$ with respect to $\hh,$} \\
3) \hbox{ $\displaystyle{\LL^0=\sum_{\a\in \Phi\setminus\{0\}}[\LL^\a,\LL^{-\a}]}$}
\end{array}
$$ and so  as above $\LL$ is completely reducible with irreducible constituents isomorphic to $\gg,$ $\u,$ $\ss$ or to the  trivial module.
 Since $\sum_{i\in I} \gg^{(i)}\op \sum_{j\in J} \u^{(j)}\op\sum_{t\in T} \mathcal{S}^{(t)}$ is a
$\gg$-submodule of $\LL,$ there is a submodule $\dd$ of $\LL$ such that
%whose weights with respect to $\hh$ contained in $\Phi$ and
 $$\LL=\sum_{i\in I} \gg^{(i)}\op \sum_{j\in J} \u^{(j)}\op\sum_{t\in T} \mathcal{S}^{(t)}\op \dd.$$ But for each nonzero $\a\in \Phi\setminus\{0\} ,$ $$\LL^\a=\mathfrak{L}^\a=\sum_{i\in I} (\fg^{(i)})^\a\op \sum_{j\in J} (\v^{(j)})^\a\op\sum_{t\in T}(\mathfrak{s}^{(t)})^\a\sub\sum_{i\in I} \gg^{(i)}\op \sum_{j\in J} \u^{(j)}\op\sum_{t\in T}\mathcal{S}^{(t)}.$$  This means that $\dd$ is a trivial $\gg$-module.
% Note that each irreducible $\gg$-module is a highest weight module.
  Now considering (\ref{final}) and  using the fact that vector spaces are flat, we may assume $$\LL=(\gg\ot \aa)\op(\ss\ot \bb)\op(\u\ot \cc)\op \dd,$$ in fact  we identify $\gg\ot \aa,$ $\ss\ot\bb$ and $\u\ot \cc$ with subspaces of $\fg\ot\aa,$ $\frak{s}\ot\bb$ and $\frak{u}\ot\cc$ respectively. Now using the same argument as in \cite[Lem. 3.6 and Pro. 3.10]{you1}, $\aa\op\bb\op\cc$ is equipped with a superalgebraic structure derived from the Lie superalgebraic structures on $\LL$ and $\fl,$ $$\fl=((\fg\ot \aa)\op(\mathfrak{u}\ot \bb)\op(\mathfrak{s}\ot \cc))+ \dd$$ and the stated properties in Theorem \ref{main} hold.

 \subsubsection{Proof of Theorem \ref{main}.} We recall that $I_0\sub I$ and $J_0\sub J$ are finite subsets  with $|I_0|=m$ and $|J_0|=n.$ Take $\Lam$ and $\Gamma$ to be index sets with a symbol  $0$ belonging to $\Lam\cap \Gamma$ such that  $\{I_\lam\mid\lam\in\Lam\}$ (resp. $\{J_\gamma\mid\gamma\in \Gamma\}$) is the set of finite subsets of $I$ (resp. $J$) containing $I_0$ (resp. $J_0$). For $(\lam,\gamma)\in\Lam\times\Gamma,$ set
{\small $$\begin{array}{ll}
\Psi^{^{\lam,\gamma}}:=\Psi\cap\hbox{span}_\bbbz\{\ep_i,\d_p\mid i\in I_\lam,p\in J_\gamma\},&
\D_1^{^{\lam,\gamma}}:=\D_\frak{u}\cap\hbox{span}_\bbbz\{\ep_i,\d_p\mid i\in I_\lam,p\in J_\gamma\},\\
R^{^{\lam,\gamma}}:=R\cap\hbox{span}_\bbbz\{\ep_i,\d_p\mid i\in I_\lam,p\in J_\gamma\},&
\D_2^{^{\lam,\gamma}}:=\D_\frak{s}\cap\hbox{span}_\bbbz\{\ep_i,\d_p\mid i\in I_\lam,p\in J_\gamma\}.
\end{array}$$}
and take
{\small
$$\begin{array}{ll}
\fl^{^{\lam,\gamma}}:=\sum_{\a\in \Psi^{^{\lam,\gamma}}}\fl^{\a}+\sum_{\a\in \Psi^{^{\lam,\gamma}}}[\fl^{\a},\fl^{-\a}],&
\frak{u}^{^{\lam,\gamma}}:=\sum_{\a\in \D_1^{^{\lam,\gamma}}}\frak{u}^{\a}+\sum_{\a\in \D_1^{^{\lam,\gamma}}}\fg^{\a}\frak{u}^{-\a},\\
\fg^{^{\lam,\gamma}}:=\sum_{\a\in R^{^{\lam,\gamma}}}\fg^{\a}+\sum_{\a\in R^{^{\lam,\gamma}}}[\fg^{\a},\fg^{-\a}],
&
\frak{s}^{^{\lam,\gamma}}:=\sum_{\a\in \D_2^{^{\lam,\gamma}}}\frak{s}^{\a}+\sum_{\a\in \D_2^{^{\lam,\gamma}}}\fg^{\a}\frak{s}^{-\a}.
\end{array}$$}
Using the result of Subsection \ref{compatibility}, we find a subsuperalgebra  $\dd$   of $\fl^{^{0,0}}$ and superspaces $\aa,$ $\bb$ and $\cc$ such that the properties stated in  Theorem \ref{main} are satisfied and
$$\fl^{^{0,0}}=(\fg^{^{0,0}}\ot \aa)\op(\mathfrak{s}^{^{0,0}}\ot \bb)\op(\mathfrak{u}^{^{0,0}}\ot \cc)\op \dd,$$
moreover, for $\lam\in\Lam$ and $\gamma\in\Gamma,$
$$\fl^{^{\lam,\gamma}}=((\fg^{^{\lam,\gamma}}\ot \aa)\op(\mathfrak{s}^{^{\lam,\gamma}}\ot \bb)\op(\mathfrak{u}^{^{\lam,\gamma}}\ot \cc))+\dd.$$
Now the result follows using the same argument as in \cite[Thm. 4.1]{you1}.\qed


\begin{thebibliography}{99}
\bibitem{AABGP} B. Allison, S. Azam, S. Berman, Y. Gao and  A. Pianzola, {\it Extended affine Lie algebras and their root systems}, Mem. Amer. Math.
Soc. 603 (1997) 1--122.


\bibitem{ABG1} B.N. Allison, G. Benkart and Y. Gao,
\newblock {\em Central extensions of Lie algebras graded by finite root systems},
\newblock  Math. Ann.  {316 } (2000),  no. 3, 499--527.


\bibitem{ABG2}
B.N. Allison, G. Benkart and Y. Gao,
\newblock {\em Lie algebras graded by the root systems $BC_r,$ $r\geq2$},
\newblock  Mem. Amer. Math. Soc.  {158 } (2002),  no. 751, x+158.

\bibitem{AKY} S. Azam, V. Khalili and  M. Yousofzadeh, {\it Extended affine root systems of type BC}, J. Lie Theory 15 (1) (2005) 145--181.
%\bibitem{AYY} S. Azam, H. Yamane and M. Yousofzadeh, {\it Reflectable bases for affine reflection systems}, J.  Algebra {371} (2012) 63--93.

\bibitem{BZ}  G. Benkart and E. Zelmanov, \newblock {\em Lie algebras graded by finite root systems and intersection matrix
algebras},
\newblock Invent. Math.  {126}  (1996),  no. 1, 1--45.




\bibitem{BS} G. Benkart and O. Smirnov, \newblock {\em Lie algebras graded by the root system $BC_1$},
\newblock J. Lie theory {13} (2003), 91--132.




\bibitem{BE1} G. Benkart and A. Elduque, {\it Lie superalgebras graded by
the root systems
$C(n), D(m,n), D(2,1;\a), F(4)$ and $ G(3),$} Canad. Math. Bull. Vol. 45 (4), (2002), 509--524.
\bibitem
%[BE2]
{BE2} G. Benkart and A. Elduque,  {\it Lie superalgebras graded by the root system $A(m, n),$} J. Lie Theory {13} (2003), 387--400.

\bibitem
%[BE3]
{BE3} G. Benkart and A. Elduque,  {\it Lie superalgebras graded by the root system $B(m,n),$ Selecta Math. (N.S.) {9} (2003), no. 3, 313--360.}









\bibitem{BM} S. Berman and R. Moody, \newblock {\em Lie algebras graded by finite root systems and the intersection matrix algebras of
 Slodowy},
\newblock Invent. Math.  {108}  (1992),  no. 2, 323--347.

\bibitem{ch-w} Sh.J. Cheng and W. Wang,  Dualities and representations of Lie superalgebras, Graduate Studies in Mathematics, 144. American Mathematical Society, Providence, RI, 2012. xviii+302 pp.


\bibitem{GN} E. Garcia and E. Neher, {\it Gelfand-Kirillov dimension and local finiteness of Jordan superpairs covered by grids and their associated Lie superalgbras},  Comm. Algebra \textbf{32} (2004), no. 6, 2149–-2175.
%\bibitem{HY}
%I. Heckenberger and  H. Yamane, {\it A generalization of Coxeter groups, root systems, and Matsumoto's theorem, Mathematishe Zeitschrift,} 259 (2008), 255--276.
%\bibitem{Hof} G. Hofmann, {\it Weyl groups with Coxeter presentation and presentation by conjugation}, J. Lie Theory 17 (2007) 337--355.
%\bibitem{Hoy} K. Hoyt, Kac-Moody superalgebras of finite growth, Ph.D. thesis, UC Berkeley, Berkeley, 2007.
\bibitem{Hum} J.E. Humphreys, {Introduction to Lie algebras and
representation theory}, Spring Verlag, New York, 1972.

%\bibitem{K} V. Kac, Infinite dimensional Lie algebras, third edition, Cambridge  University Press, 1990.
\bibitem{K1} V. Kac, {\it Lie superalgebras}, Adv. Math {26} (1977), 8--96.
\bibitem{K2} V. Kac, {\it A Sketch of Lie Superalgebra Theory}, Commun. math. Phys. {53}  (1977), 31--64.
\bibitem{LN}  O. Loos and   E. Neher, {\it Locally finite root systems}, Mem. Amer. Math. Soc. {171} (2004), no. 811, x+214.
\bibitem{LN2}  O. Loos and   E. Neher, {\it  Reflection systems and partial root systems},
Forum Math. {23} (2011), no. 2, 349--411.

\bibitem{MY} J. Morita and Y. Yoshii, {\it Locally extended affine Lie algebras}, J. Algebra {301 (1)} (2006), 59--81.

\bibitem{NS} K.H. Neeb and N. Stumme,
\newblock {\em The classification of locally finite split simple Lie algebras},
\newblock J. Reine angew. Math. {533} (2001), 25--53.
\bibitem{N1} E. Neher, {\it Extended affine Lie algebras and other generalization of affine Lie algebras- a survey}, Developments and trends in infinite-dimensional Lie theory, 53--126, Prog. Math., 228, Birkhauser Boston, Inc., Boston, MA, 2011.
\bibitem{N2} E. Neher, \newblock {\em Lie algebras graded by $3$-graded root systems and
Jordan pairs covered by grids},
\newblock Amer. J.  Math. {118} (1996), 439--491.



%\bibitem{P} I. Penkov, {\it Classically Semisimple Locally Finite Lie
%Superalgebras}, Forum Math. {16} (2004), no. 3, 431--446.
%
\bibitem{Se} G.B. Seligman, Rational methods in Lie algebras,  M. Dekker Lect. Notes in
pure and appl. math. {17}, New York, 1976.


\bibitem{serg} V. Serganova, {\it On generalizations of root systems}, comm. Algebra, {24(13)}  (1996), 4281--4299.
%\bibitem{van-de} J.W. Van de Lour, {\it A classification of Contragrediant Lie superalgebras of finite growth}, Com. in Algebra, {17(8)} (1989), 1815--1841.
%\bibitem{yos2}  Y. Yoshii, \newblock {\em Root systems extended by an abelian
%group, their Lie algebras and the classification of Lie algebras
%of type $B$},
%\newblock J. Lie theory {14 } (2004), 371--394.


\bibitem{you1} M. Yousofzadeh, {\it Structure of root graded Lie algebras}, J. Lie Theory {{22}} (2012), 397--435.

\bibitem{you4} M. Yousofzadeh, {\it Central extension of root graded Lie algebras}, Publ. Res. Inst. Math. Sci. \textbf{49} (2013), no. 4, 801-–829.
\bibitem{you2}  M. Yousofzadeh, {\it Locally finite root supersystems,}  arXiv:1309.0074.

\bibitem{you3}  M. Yousofzadeh, {\it Extended affine Lie superalgebras,} arXiv:1309.3766.

\end{thebibliography}
\end{document}